\documentclass[english]{article}
\usepackage[fontsize=11pt]{scrextend}
\usepackage[utf8]{inputenc}
\usepackage{graphicx}
\usepackage{amsmath}
\usepackage{amsfonts}
\usepackage{soul}
\usepackage{tikz}
\usepackage{upgreek}
\usetikzlibrary{matrix,fit,backgrounds}
\usepackage{enumitem}   
\usepackage{amssymb}
\usepackage{comment}
\usepackage[english]{babel}
\usepackage{hyperref} 
\usepackage[top=2cm, bottom=3cm, left=2.5cm , right=2.5cm]{geometry}
\usepackage{mathrsfs}
\usepackage{dsfont}
\usepackage{color}
\usepackage{flushend}
\usepackage{mathtools}
\usepackage{fancyhdr}
\usepackage[thmmarks,amsmath]{ntheorem}
\usepackage{stmaryrd}
\usepackage[nottoc]{tocbibind}
\usepackage[toc,page]{appendix}
\usepackage{appendix}
\usepackage{caption}
\usepackage{xcolor}
\usepackage{float}
\usepackage{textcomp}
\usepackage{bm}  
\usepackage{lipsum}
\usepackage[normalem]{ulem}
\usepackage{titlesec}
\usepackage{braket} 
\usepackage{esint} 
\usepackage{centernot}

\titleformat{\paragraph}[block]
{\normalsize\bfseries}{}{0pt}{}

\titlespacing{\paragraph}{0pt}{*1.5}{*0.5}

\renewcommand{\textbf}[1]{\begingroup\bfseries\mathversion{bold}#1\endgroup}

\newcommand{\R}{\mathbb{R}}

\renewcommand{\P}{\mathbb{P}}

\newcommand{\mcal}[1]{\mathcal{#1}}

\newcommand{\mbb}[1]{\mathbb{#1}}

\newcommand{\msf}[1]{\mathsf{#1}}
\newcommand{\mbf}[1]{\mathbf{#1}}
\newcommand{\mfrak}[1]{\mathfrak{#1}}
\newcommand{\mrm}[1]{\mathrm{#1}}
\newcommand{\tr}{\mathrm{tr}}

\newcommand{\ii}{\mathrm{i}}
\newcommand{\diff}{\mathrm{d}}

\newcommand{\tend}[1]{\underset{#1}{\longrightarrow}}

\newcommand{\defi}{\overset{(\mrm{def})}{=}}
\theoremstyle{plain}
\newtheorem{Thm}{Theorem}[section]
\newtheorem{Prop}[Thm]{Proposition}
\newtheorem{Lem}[Thm]{Lemma}
\newcommand{\vertiii}[1]{{\left\vert\kern-0.25ex\left\vert\kern-0.25ex\left\vert #1 
		\right\vert\kern-0.25ex\right\vert\kern-0.25ex\right\vert}}

{
	\theorembodyfont{\normalsize}
	\theoremsymbol{\ensuremath{\square}}
	\newtheorem*{Pro}{Proof}
	
}
\newtheorem{Def}[Thm]{Definition}
\newtheorem{Cor}[Thm]{Corollary}
{	
	\theorembodyfont{}
	\newtheorem{Rem}[Thm]{Remark}
}

{	
	\theorembodyfont{\normalsize}
	
}
\newtheorem{assu}[Thm]{Assumptions}

\begin{document}
	\title{Large deviations at the edge for 1D gases and tridiagonal random matrices at high temperature
}
	\author{Charlie Dworaczek Guera\footnote{KTH Stockholm, \newline
			\textit{email:} chadg@kth.se}\and Ronan Memin\footnote{Ecole Normale Supérieure, France \newline \textit{email:} ronan.memin@ens.fr}}
	\maketitle
	
	\begin{abstract}We consider a model for a gas of $N$ confined particles subject to a two-body repulsive interaction, namely the one-dimensional log or Riesz gas. We are interested in the so-called \textit{high temperature} regime, \textit{ie} when the inverse temperature is given by $\beta_N=2\alpha/N$ for some $\alpha>0$. We establish, in the log case, a large deviation (LD) principle and moderate deviations estimates  for the largest particle $x_\mathrm{max}$ when appropriately rescaled. Our result is in the continuity of \cite{arous2001aging,pakzad2020large} where such estimates were shown for the largest particle of the $\beta$-ensemble at fixed $\beta_N=\beta>0$ and $\beta_N\gg N^{-1}$ respectively. We show that the corresponding rate function is the same as in the case of iid particles. We also provide LD estimates in the Riesz case. 
    Additionally, we consider related models of symmetric tridiagonal random matrices with independent entries having Gaussian tails; for which we establish the LD principle for the top eigenvalue. In a certain specialization of the entries, we recover the result for the largest particle of the log-gas. We show that LD are created by a few entries taking abnormally large values.
	\end{abstract}
	
	\tableofcontents
	\section{Introduction}
	\subsection{Setting of the problem}

	
	Let $N\geq1$, $\alpha>0$, and $V:\R\to \R$ continuous. We are interested in the Gibbs probability measure given by
	\begin{equation}\label{def:mesureparticules}
		\diff\P_N(\underline{x})=p_N(\underline{x})\diff x_1\dots \diff x_N\hspace{0,5cm}\text{ with }\hspace{0,5cm}p_N(\underline{x})\defi\dfrac{1}{\mcal{Z}_N}\prod_{i<j}^{N}e^{-\frac{2\alpha}{N}g(x_i-x_j)}\cdot\prod_{i=1}^{N}e^{-V(x_i)},
	\end{equation}
	where the two-body interaction $g$ is either given by the logarithmic repulsion $g(x)=-\log|x|$ or by the Riesz repulsion $g(x)=|x|^{-s}$ with $0<s<1$. The \textit{partition function} $\mcal{Z}_N$ is the normalization constant which makes $\P_N$ a probability measure on $\R^N$, namely:
	\begin{equation}\label{def:partfunction}
		\mcal{Z}_N=\int_{\R^N}\prod_{i<j}^{N}e^{-\frac{2\alpha}{N}g(x_i-x_j)}\cdot\prod_{i=1}^{N}e^{-V(x_i)} \diff x_i.
	\end{equation}
    The measure \eqref{def:mesureparticules} is called the (one dimensional) log gas (or $\beta$-ensemble) in the case $g=-\log$, and the Riesz gas in the case $g=|.|^{-s}.$ We are interested in this paper in the regime where $N\beta_N\defi2\alpha$ is fixed (independent of $N$), also known as the \textit{high temperature regime}.
    When $V$ is assumed to grow fast enough to infinity (see the discussion after Assumptions \ref{assumptions}), $\mathcal{Z}_N$ is finite. As a consequence of \cite{Garcia} that if $\underline{x}$ is distributed with respect to $\P_N$, the \textit{empirical measure}:
	$$\mu_N\defi\dfrac{1}{N}\sum_{i=1}^N\delta_{x_i}$$
    converges almost surely towards the minimizer of a functional $\mcal{E}$ over real probability measures $\mathcal E:\mathcal{M}_1(\R)\to [0,+\infty]$. It is given by:
	\begin{equation}
		\label{def:energie}
		\mcal{E}(\mu)\defi \int_\R V(x)\diff\mu(x)+\alpha\iint_{\R^2}g(x-y)\diff\mu(x)\diff\mu(y)+\int_\R\log\dfrac{\diff\mu(x)}{\diff x}\diff\mu(x),
	\end{equation}
	for $\mu$ absolutely continuous with respect to Lebesgue measure and $\mcal{E}(\mu)\defi +\infty$ otherwise. In \eqref{def:energie}, the last term is called the \textit{entropy} of the measure $\mu$ and only differs from the physical entropy by a minus sign. Minimizing $\mcal{E}$ can then be seen as finding a compromise between a minimization of the energy (composed of the two first terms) and a maximization of the physical entropy. The latter is negligible when the inverse temperature of the Gibbs measure $\beta_N$ satisfies $\beta_N\gg N^{-1}$, but has the same order of magnitude as the energy when the correlation between particles is small/the temperature is high (see \cite[(34)]{dean2008extreme} for a discussion). 
	
	We call the minimizer of $\mcal{E}$, the \textit{equilibrium measure} and denote it $\mu_{\mrm{eq}}$ (we omit the $V$ and $\alpha$-dependence since these parameters remain fixed throughout the entire article). It is absolutely continuous (with respect to Lebesgue measure) with p.d.f. $\rho_{\mrm{eq}}$ satisfying:
	\begin{equation}\label{eqmeasure:charac}V(x)+2\alpha \int_\R g(x-y)\rho_{\mrm{eq}}(y)dy+\log\rho_{\mrm{eq}}(x)=\lambda_{\mrm{eq}}\hspace{1cm}x-ae
	\end{equation}
	where $\lambda_{\mrm{eq}}$ is a constant (see \cite[Proposition A.2]{lambert2021poisson}). In this context, $\rho_{\mrm{eq}}$ is supported on $\R$ because of the presence of the entropy in the minimizing equation. This is a major difference with the case where $\beta$ is constant, where the equilibrium measure is compactly supported. This equation can be rewritten as
	\begin{equation}\label{eqmeasure:exp}
		\forall x\in\R,\quad \rho_{\mrm{eq}}(x)=\exp\Big(-V(x)-2\alpha U^{\rho_{\mrm{eq}}}(x)+\lambda_{\mrm{eq}}\Big),\hspace{1cm}U^{\rho_{\mrm{eq}}}(x)\defi \displaystyle\int_\R g(x-y)\rho_{\mrm{eq}}(y)\diff y.
	\end{equation}
	The function $U^{\rho_{\mrm{eq}}}$ called the Riesz or logarithmic (according to the choice of $g$) potential, has recently been identified as the Markov-Krein transform of simple distributions for some specific examples of $V$ and $g=-\log|.|$ in \cite{nakano2024classical}. 
    The density $\rho_{\mrm{eq}}$ also admits explicit expressions, such as in the case $V_G(x)=x^2/2$ and $g=-\log|.|$ where it is known as the \textit{Askey-Wimp-Kerov} distribution (after \cite{askey1984associated} and \cite[§ 8]{kerov1998interlacing}). In the notations of \cite[Eq. (18)]{AllezBouchaudGuionnet}, it reads:
	$$\rho_{V_G}(x)=\dfrac{1}{\sqrt{2\pi}\Gamma(1+\alpha)}\dfrac{1}{|D_{-\alpha}(\ii x)|^2},\hspace{2cm}D_{-\alpha}(z)\defi  \frac{1}{\Gamma(\alpha)}e^{-\frac{z^2}{4}}\displaystyle\int_0^{+\infty}t^{\alpha-1}e^{-\frac{t^2}{2}-zt}\diff t.$$
	Because the limiting empirical measure $\mu_{\mrm{eq}}$ is supported on the whole real line, the largest particle 
    $$x_{\max}\defi  \displaystyle\max_{1\leq i\leq N} x_i$$ goes to $+\infty$ almost surely. More precisely, \cite[Theorems 1.7, 1.8]{lambert2021poisson} establish (in the general context of log and Riesz interactions in the high temperature regime) that for $V=|.|^ d$, $x_{\max}$ is almost surely of order $(\log N)^{1/ d}$, by showing that the local statistics at the edge are given by a Poisson point process. Considering the particular case $\beta_N=0$ of independent particles, these results are well known, see for example \cite{vivo2015large}. This shows that the correlations introduced by the interaction $g$ in \eqref{def:mesureparticules} do not affect the nature of the statistics of the largest particle in this regime.

    In this paper, we show a LDP for the largest particle $x_{\max}$, namely a quantification of the probability of a deviation of $x_{\max}$ from its typical behavior in the large $N$ limit.\\

    In the $\log$ gas case, \textit{ie} when $g=-\log$, another approach to analyze the measure \eqref{def:mesureparticules} is given by its link with random matrices. Indeed, specifying $V(x)=x^2/2$, it was shown in \cite[Theorem 2.12]{dued} that $\P_N$ is the joint law of (unordered) eigenvalues of the tridiagonal matrix:
    \begin{equation}\label{eq:dued}
        T= \begin{pmatrix}
    a_1 & b_1 &         &       &\\
    b_1 & a_2 & b_2     &             \\
        & b_2 & \ddots  & \ddots      \\
        &     & \ddots  &  \ddots       & b_{N-1}      \\
    &     &         &b_{N-1}  & a_N
\end{pmatrix},
    \end{equation}
    where $(a_i)_{1\leq i \leq N}$, $(b_i)_{1\leq i \leq N-1}$ are independent families of independent random variables, with $a_i \sim \mathcal{N}(0,1)$ and $b_i \sim \chi_{\beta_N(N-i)}/\sqrt{2}$, where for $a >0$ the law $\chi_a$ has density \begin{equation}
    \label{def:chi}
        f_a(x) = \frac{2^{1-a/2}}{\Gamma(a/2)}x^{a-1}e^{-x^2/2}\mathbf{1}_{x>0}.
    \end{equation}
    With the representation \eqref{eq:dued}, it is natural to ask which rare events involving the entries of $T$ are responsible for a large deviation of the top eigenvalue. We establish in the sequel that such large deviations are driven by a small number of neighbouring entries taking large values.
    \subsection{Notation}
   For $(u_N),(v_N)$ real sequences with $v_N\neq 0$, we write $u_N= o(v_N)$ (resp. $u_N=O(v_N)$, resp. $u_N\sim v_N$) if $u_N/v_N\rightarrow0$ (resp. is bounded, $u_N=v_N+o(v_N)$). We use the same notations for functions. We write $\mathcal{M}_1(E)$ for the set of Borel probability measures on $E$. For a finite set $X$, we denote by $\mrm{Card}(X)$ its cardinal. For $a<b$ integers we denote $\llbracket a,b\rrbracket=\{a,a+1,\ldots,b\}$, and $a\wedge b =\min(a,b)$.
    
    We denote by $\fint_\R f(x)\diff x$, the Cauchy principal-value of $f$.
    To distinguish the cases in several proofs, we use the notation $\vartheta = \log(1+|\cdot|)\mathbf{1}_{g=-\log}$
	\subsection{Main results}
	\begin{assu}\label{assumptions} The potential $V$ satisfies:
	$V(x)=\theta|x|^ d + \varphi(x)$ where $d\in [1,+\infty)$, $\theta>0$ and $\varphi(x)\mathop{=}\limits_{|x|\to\infty}o(x^{ d})$, where $\varphi$ is $\mathcal{C}^2$ outside of a compact set, with $\varphi'(x)\mathop{=}\limits_{|x|\to\infty}o(x^{ d-1})$ and for any $a<b$, $$\sup_{y\in [a,b]}\varphi''(x+y)\mathop{=}\limits_{|x|\to\infty}o(x^{ d-2}).$$
	\end{assu}
    The assumptions above are quite general. The bound on $\varphi''$ allows us to use the results from \cite{lambert2021poisson} about the convergence of the edge towards a Poisson point process when $N\rightarrow\infty$, see Lemma \ref{lem:lambertedge}.
    Assumptions \ref{assumptions} also imply that there exist $\delta>0$ and $c\in\R$ such that for all $x\in\R$, $$V(x)-(2\alpha +1+\delta)\log(1+|x|)\geq c,$$ ensuring that the probability measure \eqref{def:mesureparticules} is well defined.
    
	Define $I_ d:\R\mapsto[0,+\infty]$ by
	\begin{equation}
    \label{def:ratefunction}
		I_ d(x)\defi\begin{cases}
			x^ d
			-1\hspace{1cm}&\mrm{if}\;x\geq1
			\\+\infty\;&\mrm{if}\;x<1
		\end{cases}.
	\end{equation}
	This function is a good rate function, meaning that it is lower semi-continuous and its level sets are compact.
    Our main result establishes a full LDP for the largest particle of the $\log$ gas. In what follows, we set:
    $$\overline{x_\mrm{max}}\defi \dfrac{x_\mrm{max}}{(\theta^{-1}\log N)^{1/ d}}$$
	\begin{Thm}[LDP for $\beta$-ensembles]
		\label{thm:LDP}
		Under Assumptions \ref{assumptions}, let $g=-\log$. The law of
        $\overline{x_\mrm{max}}$ under $\mbb{P}_N$ satisfies a LDP at speed $\log N$ with good rate function $I_ d$, \textit{ie} the following estimates hold:
		\begin{itemize}
			\item For any open set $O\subset \R$, 
			\begin{equation}
				\label{eq:LDPLower}
				-\inf_{x\in O}I_ d(x)\leq\liminf_{N\rightarrow\infty}\dfrac{1}{\log N}\log\mbb{P}_N\left(\overline{x_\mrm{max}}\in O\right)
			\end{equation}
			\item For any closed set $F\subset \R$, 
			\begin{equation}
				\label{eq:LDPUpper}
				\limsup_{N\rightarrow\infty}\dfrac{1}{\log N}\log\mbb{P}_N\left(\overline{x_\mrm{max}}\in F\right)\leq-\inf_{x\in F}I_ d(x).
			\end{equation}
		\end{itemize}
	\end{Thm}

    In the Riesz case, we prove the following LD estimates.

\begin{Thm}[Riesz case]\label{thm:rieszLDP}Under Assumptions \ref{assumptions}, let $g= |.|^{-s}$. Then,
	\label{thm:Right deviation}
		for all $x>1$,
			\begin{equation*}
				\limsup_{N\rightarrow\infty}\dfrac{1}{\log N}\log\mbb{P}_N\left(\overline{x_\mrm{max}}>x\right)\leq-I_ d(x),
			\end{equation*}
            and
			\begin{equation*}
				\lim_{\delta\to 0}\liminf_{N\rightarrow\infty}\dfrac{1}{\log N}\log\mbb{P}_N\Big(\overline{x_\mrm{max}}\in (x-\delta,x+\delta)\Big)\geq-I_ d(x).
			\end{equation*}
            For all $x< \left(\dfrac{1-s}{1+s}\right)^{1/ d}$,
            \begin{equation}\label{eq:leftriesz}
            \lim_{N\rightarrow\infty}\dfrac{1}{\log N}\log\mbb{P}_N\left(\overline{x_\mrm{max}}<x\right)=I_ d(x)=-\infty.
			\end{equation}
	\end{Thm}
    
	\begin{Rem}
		\begin{itemize}
        \item  Although we expect the full LDP for the largest particle in the Riesz case to hold, Theorem \ref{thm:rieszLDP} says nothing about the rate of decay of $\mbb{P}_N\left(\overline{x_\mrm{max}}<x\right)$ in the region $\left(\dfrac{1-s}{1+s}\right)^{1/ d}$ $\left(\dfrac{1-s}{1+s}\right)^{-1/ d}\leq x<1$. To do so, extending \eqref{eq:leftriesz} to all $x<1$ would be enough. 
			\item Our result only covers the case of potentials with power growth. However, it can be seen that for potentials growing faster to infinity, the LDP looks different. For example, in the case of iid particles $(\alpha=0)$ and potential $V=\cosh$, the largest particle satisfies a LDP at speed $\log N$ with rate function given by
			$$I_{\cosh}(x)\defi\begin{cases}
				&+\infty\hspace{1cm}\text{if }x<1
				\\&0\hspace{1,5cm}\text{if }x=1
				\\&+\infty\hspace{1cm}\text{if }x>1
			\end{cases}$$
            The explanation for this rate function to be trivial is because the decay speed of the log probability of a deviation depends on the deviation. For example, one can see that for $x>0$,
            $$\dfrac{1}{\log N}\log\mbb{P}_N\left(\dfrac{x_\mrm{max}}{\log(2\log N)}>1+x\right)=-(\log N)^{x}[1+o_{N}(1)].$$
			\item The analogue of this result for $\beta$-ensembles for $\beta>0$ are often used in truncation arguments such as in \cite[Proposition 2.2]{borot2013asymptotic} and could be applied to extend the asymptotic expansion of $\log\mcal{Z}_N$ proven in \cite{Dwo} to more general potentials.
            \item The previous results are enough to justify the LD for the law of $\overline{x_\mrm{min}}\defi\dfrac{x_\mrm{min}}{(\theta^{-1}\log N)^{1/ d}}$. Indeed, the law of $\overline{x_\mrm{min}}$ is equal to the law of $-\overline{x_\mrm{max}}$ under the choice of potential $W(x)=V(-x)$.
            \end{itemize}
	\end{Rem}


    Establishing Theorems \ref{thm:LDP} and \ref{thm:rieszLDP}, we derive a refined upper bound in the case of a left deviation for $x_{\max}$.

    Contrary to a right deviation, a left one (\textit{i.e.} the case where the largest particle is smaller from its typical position) imposes a constraint on a mesoscopic number (of the form $N^\varepsilon$ with $0<\varepsilon<1$) of particles, making this kind of deviations much more unlikely than the right ones.
	\begin{Prop}[Refined upper bound]\label{prop:preciseupper}
		Under Assumptions \ref{assumptions}, set $x_0 = 0$ for $g=-\log$, and $x_0=2s/(1+s)$ for $g=|.|^{-s}$. Then for all $x_0<x<1$,
		\begin{equation*}
			\limsup_{N\rightarrow\infty}\dfrac{1}{\log N}\log\log \Big(\mbb{P}_N\left(\overline{x_\mrm{max}}<1-x\right)^{-1}\Big)\geq 1-(1-x)^ d.
		\end{equation*}
	\end{Prop}

    \begin{Rem}
        We expect this bound to remain true for $0<x<x_0$ in the Riesz case as the edge should behave the same as the edge of $N$ iid particles. On the event $\{\overline{x_{\max}}\leq 1-x\}$, we should observe the displacement of a mesoscopic number of particles of order $N^{\varepsilon}$ with $\varepsilon=1-(1-x)^ d$. For this reason, we expect the lower bound to hold too. However, since we were only interested in showing that $I_d$ was infinite on $(-\infty,1)$ we didn't try to investigate it further.
    \end{Rem}

    In this paragraph, take $\theta=1$ for simplicity. As a consequence of \cite[Theorems 1.7, 1.8]{lambert2021poisson}, Lambert establishes that for some sequence $(E_N)_{N\geq 1}$ satisfying $E_N\sim (\log N)^{1/ d}$ (see Lemma \ref{lem:lambertedge}), $ d E_N^ d\big(x_\mathrm{max}/E_N-1\big)$ converges in law towards a Gumbel random variable. In particular, $E_N\sim (\log N)^{1/ d}$ is identified as the typical edge of a configuration, and for $x>0$, the probability $\P_N\big(|x_\mathrm{max}-(\log N)^{1/ d}|>x/(\log N)^\gamma\big)$ has to go to zero if $\gamma<1- d^{-1}$. This regime is called the \textit{moderate deviations} regime and the quantification of these probabilities is the object of the following proposition.

    \begin{Prop}[Moderate deviations]\label{thm:moderate}
			Under Assumptions \ref{assumptions}, assuming in addition that $\varphi(x)=O_{|x|\to \infty}(1)$, for all $- d^{-1}<\gamma<1- d^{-1}$ and $x>0$ the following estimates hold:
            \begin{equation*}
                \limsup_{N\rightarrow\infty}\dfrac{1}{(\log N)^{1-\gamma- d^{-1}}}\log\mbb{P}_N\Big(x_{\mrm{max}}\geq (\theta^{-1}\log N)^{1/ d}+\dfrac{x}{(\log N)^\gamma}\Big)\leq -xd\,,
            \end{equation*}
            and for $g=-\log$, for all $x>0$:
        \begin{equation*}
			\limsup_{N\rightarrow\infty}\dfrac{1}{(\log N)^{1-\gamma- d^{-1}}}\log\log\Big(\mbb{P}_N\Big(x_{\mrm{max}}\leq (\theta^{-1}\log N)^{1/ d}-\dfrac{x}{(\log N)^\gamma}\Big)^{-1}\Big)\geq xd\,.
		\end{equation*}
	\end{Prop}
    \begin{Rem}
      We didn't try to investigate the lower bounds for the moderate deviations but we expect our lower bound to be optimal.
    \end{Rem}
	    In Theorem \ref{thm:LDPtridiag_intro}, we establish the following LDP for (periodic or not) tridiagonal matrices with independent entries having Gaussian tails. Let $(a_1,\ldots,a_N),(b_1,\ldots,b_N)\in \R^N$, we consider the symmetric tridiagonal matrix:
    \begin{equation}
        \label{def:T}
        T= \begin{pmatrix}
    a_1 & b_1 &         &       &\\
    b_1 & a_2 & b_2     &             \\
        & b_2 & \ddots  & \ddots      \\
        &     & \ddots  &  \ddots       & b_{N-1}      \\
    &     &         &b_{N-1}  & a_N
\end{pmatrix}
    \end{equation} 
    and the corresponding symmetric \textit{periodic} tridiagonal matrix 
\begin{equation}
        \label{def:T_per}
        T_\mathsf{per}= \begin{pmatrix}
    a_1 & b_1 &         &       & b_N \\
    b_1 & a_2 & b_2     &             \\
        & b_2 & \ddots  & \ddots      \\
        &     & \ddots  &  \ddots       & b_{N-1}      \\
    b_N &     &         &b_{N-1}  & a_N
\end{pmatrix}.
    \end{equation} 
    We assume the entries of these matrices to have \textit{uniform Gaussian tail}:
\begin{Def}[Uniform Gaussian tail]\label{def:gaussian_tail}
    We say that the family of real random variables $(X_i)_{i\in I}$, where $I$ is any index set, has uniform Gaussian tail if for all $i\in I$,
    $$ -\log\P\left(|X_i|>x\right) \mathop{=}\limits_{x\to +\infty} x^2/2 +o(x^2),\quad \text{and} \quad -\log\P\left(X_i>x\right)\mathop{=}\limits_{x\to +\infty}x^2/2+\mathop{o}(x^2),$$
    where the reminder term $o(x^2)$ is bounded independently of $i\in I$.
    \end{Def}

We show the following.

\begin{Thm}\label{thm:LDPtridiag_intro}
        Let $(a_i^{(N)})_{i\geq 1},(b_i^{(N)})_{i\geq 1}$ be independent sequences of independent variables, possibly depending on $N$. Assume that for $N$ large enough $(a_i^{(N)})_{i\geq 1}$ and $(\sqrt{2}b_i^{(N)})_{i\geq 1}$ have uniform Gaussian tail.
        
        Consider the tridiagonal or periodic tridiagonal model $$J=T(a_1^{(N)},\ldots,a_N^{(N)},b_1^{(N)},\ldots,b_{N-1}^{(N)})\quad\text{or}\quad J=T_\mathsf{per}(a_1^{(N)},\ldots,a_N^{(N)},b_1^{(N)},\ldots,b_{N}^{(N)}).$$ Then, $\lambda_\mathrm{max}(J)/\sqrt{2\log N}$ satisfies a LDP with rate function $I_2$ given by \eqref{def:ratefunction}.
    \end{Thm}
    \begin{Rem}\label{rem:few_entries}
    \begin{itemize}
        \item  Taking for $i\geq 1$, $X_i \sim \chi_{u_i}$, with $0 < u_i < C$ for some $C$ independent of $i$, the sequence $(X_i)$ satisfies the previous conditions.
        \item As mentioned earlier, considering the model $T$ with $(a_i)_{1\leq i \leq N}$ iid standard Gaussians, and $\sqrt{2}b_i^{(N)}$ distributed with respect to $\chi_{\beta_N(N-i)}$, by \cite[Theorem 2.12]{dued}, the law of the spectrum of $T$ is distributed with respect to the high temperature Gaussian $\beta$-ensemble,
        and one recovers Theorem \ref{thm:LDP} in this special case. 
        \item The study of the periodic model $T_\mathsf{per}$ is motivated by the theory of integrable systems, as detailed in the Subsection \ref{sec:literature} \textit{Sparse and tridiagonal matrices}. Indeed  taking $(a_i)_{1\leq i \leq N}$ iid standard Gaussians and $(\sqrt{2}b_i)_{1\leq i\leq N}$ iid distributed with respect to $\chi_{2\alpha}$, $T_\msf{per}$ arises as the Lax matrix of the so-called periodic Toda chain in thermal equilibrium.
    \end{itemize}
    \end{Rem}

	\subsection{Comparison with the literature}\label{sec:literature}
	
	\paragraph{Links with extremal value of random matrices}
	
	The first LD result in random matrix theory was established by Ben-Arous and Guionnet \cite{arous1997large}, where the authors showed that the law of the empirical measure of the fixed temperature $\beta$-ensemble --- which includes the classical ensembles GOE, GUE and GSE --- satisfies a LDP. They provided an expression for the rate function, relying on the explicit joint law of eigenvalues.
	
 	A LDP for the largest eigenvalue was first shown for the GOE in \cite{arous2001aging} in the context of the study of the aging phenomenon in spin glass models. This result was then extended in \cite[Theorem 2.6.6]{anderson2010introduction} for general $\beta$-ensembles (at fixed $\beta>0$) whose partition function satisfies a certain rate condition. This type of LDP was later shown, proving a continuity result for the spherical integrals, for rank-one perturbations of GOE/GUE in \cite{maida2007large} and then extended for the joint law of the extremal eigenvalues in \cite{benaych2012large}. The case of critical $\beta$-ensembles was also studied in \cite{fan2015convergence}. The LDP for the largest eigenvalue was recently further generalized in \cite{guionnet2024asymptoticexpansionpartitionfunction}, where the authors consider $\beta$-ensembles on a complex contour homeomorphic to the real line. LDP for the law of the maximal eigenvalue of Wigner matrices with variances profiles were shown in \cite{husson2022large,ducatez2024large}. In \cite{BordCaputo}, the authors establish a LDP for the empirical measure for a Hermitian matrix $X_N$ whose entries have stretched exponential tails vanishing more slowly than Gaussians, relying on the fact that LD are caused by a small number of entries being abnormally large.  In \cite{Augeri}, Augeri built on the strategy this idea to deduce a LDP for the largest eigenvalue for the same matrix model. In \cite{guionnet2020large}, the authors proved a LDP for the largest eigenvalue in the case of general sharp sub-Gaussian entries. The sharpness assumption was then relaxed in \cite{augeri2021large} where the authors considered entries with symmetric distributions. This symmetric condition was then relaxed in \cite{cook2023full}.

	\paragraph{The $\beta$-ensembles at high temperature}
	When the interaction in \eqref{def:mesureparticules} is given by $g=-\log$, a lot is known in the regime $\beta_N\rightarrow0$. A pioneer work at high temperature was \cite{AllezBouchaudGuionnet} where the authors considered the Gaussian $\beta$-ensembles with $N\beta=2\alpha $ and the so-called Gauss-Wigner crossover when $\alpha\rightarrow 0$ and $\alpha\rightarrow\infty$. Interpolations between the Tracy-Widom and Gumbel distributions, which govern respectively the fluctuations at the edge for strongly correlated and weakly correlated systems, have been studied in \cite{johansson2007gumbel,allez2014tracy}. In \cite{allez2014sine} the authors showed that the Sine$_\beta$ point process, which describes the limiting microscopic behavior of $\beta$-ensembles for fixed $\beta>0$, converges as $\beta\to 0$ to a Poisson point process. It was later shown in \cite{benaych2015poisson} that the local satistics converge, in the case of a quadratic potential, towards those of a Poisson point process similarly to the case of iid particles ($\alpha=0$). This result was then extended to the edge, general potential, general interaction and general configuration spaces in \cite{pakzad2018poisson,NakanoTrinh20,lambert2021poisson,padilla2024emergence,padilla2025poisson}. The LDP for the empirical measure for general gases (including the log and Riesz gases considered in this paper) at high-temperature was proven in \cite{Garcia}. LDP for the empirical measure were shown for Jacobi $\beta$-ensembles (resp. Laguerre) in \cite{ma2023limit} (resp. \cite{ma2023unified}) for $\beta_N\gg \log N/N$; and for the extremal eigenvalues in the same temperature regime in the Jacobi case \cite{lei2023large}. In this regime for $\beta_N$, the entropy term does not appear in the rate function for the empirical measure and one recovers the rate function in the fixed $\beta>0$ case. Central limit theorems in the regime $N\beta\to 2\alpha >0$ were established in \cite{NakanoTrinh,hardy2021clt, DwoMemin} and the asymptotics of the partition function was established in \cite{Dwo}.
	
	\paragraph{Sparse and tridiagonal random matrices}
	Many random graph models can be studied through their adjacency matrices, for example regular graphs or Erdős–Rényi graphs. A lot is known about the extremal eigenvalues of these matrices \cite{benaych2019largest,benaych2020spectral,altducatezknowles,alt2023poisson,augeri2023large}. Upper tail LD bounds were shown for the two largest eigenvalues of sparse  Erdős–Rényi ($p\rightarrow0$) in \cite{bhattacharya2020upper,cook2020large} in the regime $n^{-1/2}\ll p\ll1 $. The sparsity regime was then extended to $p$ such that $\log \tfrac{1}{np} \ll\log n$ and $np\ll \sqrt{\log n/\log\log n}$ in \cite{bhattacharya2021spectral}. Upper tail large deviations for the intermediate regime, $ \sqrt{\log n/\log\log n}\ll np\ll n^{1/2+o(1)}$ were established in \cite{basak2023upper}. In \cite{ganguly2022large,ganguly2024spectral}, the authors were able to show the LDP for a sparsity parameter $p=O(n^{-1})$ in the case of respectively Gaussian weights and general weights on edges. Although the matrix models and the techniques differ a lot, this LDP is very similar to the one obtained in this paper. The situation described in Remark \ref{rem:few_entries}, \textit{i.e.} that the LD of $\lambda_\mathrm{max}$ is produced by a bounded number of entries, is reminiscent of the conclusion of \cite{Augeri}, where the matrix of interest is a Wigner matrix whose entries have tails heavier than Gaussians. In the context of sparse graphs, in \cite{augeri2025large}, the author establishes a LDP for the empirical measure of Erdős–Rényi adjacency matrices in the regime $\log N/N\ll p\ll 1$. In \cite{han2025deviation}, LD estimates for the top eigenvalue of random Schrödinger operators and tridiagonal models related to the $\beta$-ensembles are derived. Furthermore, an application of Theorem \ref{thm:LDPtridiag_intro} is taking $(a_i)_{i\geq 1}$ iid standard Gaussians, and $(\sqrt{2}b_i)_{i\geq 1}$ iid with law $\chi_{2\alpha }$ for some $P>0$. Then, the periodic matrix $T_\mathsf{per}$, given by \eqref{def:T_per}, is distributed as the so-called Lax matrix of the periodic Toda chain (a classical instance of an integrable system), with pressure $\alpha=P>0$ at thermal equilibirum, see \cite{spohn2020generalized}. This model has attracted attention lately, via its links with the  $\beta$-ensemble at high temperature, see the works \cite{mazzuca2022mean,GMToda,mazzuca2024clt,grava2026large} investigating the links between the statistical properties of the Toda chain and the $\beta$-ensembles at high temperature.
    
	\subsection{Outline of the proofs}

    \subsubsection{Preliminaries on large deviations}
    Let us recall standard facts and definitions here. Let $(X_N)$ be a sequence of real random variables with having laws $P_N\in \mathcal{M}_1(\R)$. Let $I:\R\to[0,+\infty]$ be lower semicontinuous, meaning that for each $a\in [0,+\infty)$, $I^{-1}([0,a])\subset \R$ is closed. Let $(a_N)$ be a sequence of positive reals, with $\lim_{N\rightarrow\infty} a_N= +\infty$.
    \begin{Def}[Full/weak LDP]
    The sequence $(X_N)_{N\geq 1}$ satisfies a 
    \begin{itemize}
        \item full LDP with rate function $I$ and speed $a_N$ if for all $U\subset \R$ open,
    \begin{equation}\label{def:borneinf}
            -\inf_A I \leq \liminf_N\frac{1}{a_N}\log P_N(U),
        \end{equation}
        and for all $F\subset \R$ closed,
        \begin{equation}\label{def:bornesup}
            \limsup_N\frac{1}{a_N}\log P_N(F) \leq -\inf_A I.
        \end{equation}
        
        \item weak LDP with rate function $I$ and speed $a_N$ if \eqref{def:borneinf} holds for all $U$ open, and \eqref{def:bornesup} holds for $F\subset \R$ compact.
    \end{itemize}
    \end{Def}
    A full LDP implies a weak LDP, and the converse holds if $(X_N)_{N\geq 1}$ is exponentially tight, defined as follows.
    \begin{Def}[Exponential tightness]
        $(X_N)$ is said to be exponentially tight at speed $a_N$ if for all $M>0$, there exists a compact $K_M\subset \R$ such that
        $$ \limsup_N \frac{1}{a_N}\log P_N(K_M^\mathsf{c} ) < -M.$$
    \end{Def}
    
    This amounts to say that up to negligible error with respect to the scale of large deviations, the sequence $(X_N)_{N\geq 1}$ lives in a compact. Then, one has by Lemma \cite[Lemma 1.2.18]{dembo2009large}:
    \begin{Lem}
        If $(X_N)_{N\geq 1}$ satisfies a weak LDP with rate function $I$ and is exponentially tight at speed $a_N$, then $(X_N)_{N\geq 1}$ satisfies a full LDP with rate function $I$.
    \end{Lem}
    A way to establish a weak LDP is through the following result that can be found in \cite[Theorem 4.1.11]{dembo2009large}.
    \begin{Lem}
        $(X_N)$ satisfies the weak LDP at speed $a_N$ and rate function $I$ if for all $x\in \R$, 
        \begin{equation*}
\lim_{\delta\downarrow 0}\liminf_{N\rightarrow\infty} \frac{1}{a_N}\log P_N\big((x-\delta,x+\delta)\big)=
			\lim_{\delta\downarrow 0}\limsup_{N\rightarrow\infty} \frac{1}{a_N}\log P_N\big((x-\delta,x+\delta)\big)=-I(x).
		\end{equation*}
    \end{Lem}

    \subsubsection{Outline of the proof of Theorems \ref{thm:LDP},\ref{prop:preciseupper}}

	
		
	
 The outline of the proof of the LDP for the log case, as well as LD estimates for the Riesz gas are established in:
	\begin{itemize}
		\item[1-] Propositions \ref{prop:upperboundx>1} and \ref{prop:LUP} show respectively the large and moderate deviations upper bounds of Theorems \ref{thm:LDP}, \ref{thm:rieszLDP} and Proposition \ref{thm:moderate} for closed sets of the form $[x,+\infty)$ for $x\geq 1$, and $(-\infty, x]$ for $x<1$. The main ingredient for the upper bound for the right deviations, is an explicit bound on the marginal of the particle $x_1$, called a Wegner estimate and proven in \cite{lambert2021poisson}. For the upper bound on the left deviations, we obtain, using energy estimates as in \cite{MMS}, a pointwise bound on $p_N$ defined in \eqref{def:mesureparticules}: we show that up to a negligible remainder $\mfrak{R}_N$, it holds that for all $\underline{x}\in\R^N$: $$p_N(\underline{x})\leq e^{\mfrak{R}_N}\prod_{i=1}^N\rho_{\mrm{eq}}(x_i).$$
        In doing so, we establish on the way Proposition \ref{prop:preciseupper}.
		\item[2-] As a direct consequence of  Propositions \ref{prop:upperboundx>1}, we deduce the exponential tightness in Corollary \ref{cor:exptight} by taking $K_M\defi[-M,M]$.
		\item[3-] We establish in Proposition \ref{pro:lowerbound} the weak LD lower bound. We restrict the configuration to a deviation of only one particle while the others remain in a very likely set. Using the convergence of the local statistics towards a Poisson point process, we are able to prove that the probability of this unique particle deviating can be bounded by below as in Theorem \ref{thm:LDP}, \ref{thm:rieszLDP}.
	\end{itemize}

\subsubsection{Outline of the proof of Theorem \ref{thm:LDPtridiag_intro}}
The proof of Theorem \ref{thm:LDPtridiag_intro} relies on a truncation procedure on the entries of the matrix, similar to the procedure of \cite{Augeri}, showing that only the big entries of the matrix at stake matter for the LD. Let us describe here the case of the (nonperiodic) tridiagonal model $T_N$.

Take $x\geq 1$. The proof of the LD lower bound for $\P(\overline\lambda_{\max}\geq x)$ with $\overline\lambda_{\max}=\lambda_\mathrm{max}(T_N)/\sqrt{2\log N}$, is a consequence of the simple inequality $\max_{1\leq i \leq N}a_i^{(N)} \leq \lambda_{\max},$
see Lemma \ref{lem:inclusion_bourrine} and Remark \ref{rem:matricelower}. For the LD upper bound, we proceed as follows. For $\varepsilon>0$, we let $T_N^\varepsilon=(T_{i,j}\mathbf{1}_{ \{|J_{i,j}|\geq \varepsilon\sqrt{2\log N}\} })_{1\leq i,j\leq N}$ be the truncated matrix with entries larger than $\varepsilon\sqrt{2\log N}$. Then Lemma \ref{lem:exp_approx} shows that the normalized top eigenvalue of $T^\varepsilon_N$ is \textit{exponentially equivalent} to the one of $T_N$: more precisely, let $\overline\lambda_{\max}=\lambda_\mathrm{max}(T_N)/\sqrt{2\log N}$ and define $\overline\lambda_{\max}^{(\varepsilon)}$ similarly. Then, for any $\varepsilon>0$, we have
\begin{equation}
    \label{intro exp approx}
    \lim_{\delta\to 0}\lim_{N\rightarrow\infty}\frac{1}{\log N} \log \P(|\overline\lambda_{\max}-\overline\lambda_{\max}^{(\varepsilon)}|>\delta) = -\infty.
\end{equation}
This allows us to reduce to showing the LD upper bound for $T_N^\varepsilon$.
Now, because we assume the entries of $T_N$ to have Gaussian tails, a lot of entries in $T_N^\varepsilon$ vanish. In particular, $T_N^{\varepsilon}$ can be assumed to be a block diagonal matrix, $T^\varepsilon_N=\mathrm{diag}(B_1,\ldots,B_{k_N})$, where each block can be assumed -- at the scale of LD -- to be of finite size (\textit{ie} independent from $N$). The top eigenvalue of $T_N^\varepsilon$ is then 
\begin{equation}
    \label{intro block decomp}
    \lambda_\mathrm{\max}^{(\varepsilon)} = \max_{1\leq i \leq k_N} \lambda_\mathrm{max}(B_i).
\end{equation}
We then bound by above the largest eigenvalue of a block $B_i$ by its Euclidean norm $\|B_i\|\defi(\tr B_i^2)^{1/2}$, and conclude using Lemma \ref{lem:gaussian_tail}, stating that if $X=(X_1,\ldots,X_d)$, with $X_i$ independent with Gaussian tail, $\|X\|$ also has a Gaussian tail. Putting these arguments together, this yields the desired LD upper bound.

\paragraph{Outline of the paper.}
    

     In Section \ref{sec:logriesz}, we consider the $\log$ and Riesz gases, and establish a LDP in the first case, see Theorem \ref{thm:LDP} and LD estimates for the second in Theorem \ref{thm:rieszLDP}. We also provide further estimates (including moderate deviations) in Theorems \ref{prop:preciseupper} and \ref{thm:moderate}.
     
    In Section \ref{sec:matrix},
    we establish a LDP for the top eigenvalue of symmetric tridiagonal matrices having independent entries, and Gaussian tails, see Theorem \ref{thm:LDPtridiag_intro}. Because of the tridiagonal matrix representation \eqref{eq:dued} of the $\log$ gas, we recover in particular the result of Theorem \ref{thm:LDP} in the case of a quadratic potential, and give an interpretation of the events producing a LD of the top eigenvalue.
	
\subsection*{Acknowledgements}
The authors would like to thank the organizers of the “Random Matrices and Free probability” summer school organized in Toulouse in June 2024 where this project was initiated. The authors also thank Gaultier Lambert, Michel Pain and Jonathan Husson for fruitful discussions and suggestions. We thank two anonymous referees for fruitful remarks and helping us improve the present paper. \\
C.D.G. acknowledges the support of the starting grant 2022-04882 from the Swedish Research Council and of the starting grant from the Ragnar Söderbergs Foundation.\\
R.M. has received support from the LabEx CIMI - ANR-11-LABX-0040, and from the Major Research Program of PSL Research University "Statistical Physics and Mathematics" launched by PSL Research University and implemented by ANR-10-IDEX-0001.
    
\section{Proof of Theorems \ref{thm:LDP}, \ref{thm:rieszLDP}, Propositions \ref{prop:preciseupper}, \ref{thm:moderate}.}\label{sec:logriesz}
We first notice that in both $\log$ gas and Riesz cases, we can assume $\theta=1$ in Assumtions \ref{assumptions}, \textit{ie} $V(x)=|x|^ d + \varphi(x)$, via considering the system $\underline{y}:=\theta^{1/ d}\underline{x}$, which amounts to replacing $\varphi$ by $\varphi(\cdot/\theta^{1/ d})$, and $g$ by $g(\cdot/\theta^{1/ d})$, for which the proofs are unchanged. In this section, we therefore assume $\theta=1$.

	\subsection{Large deviation upper bounds, exponential tightness}
	
	\subsubsection{Upper bound: right-deviation}
	We start the proof by showing the upper bound in the case $x\geq1$ (\textit{i.e.} when the rate function $I_ d(x)$ is finite) and in the moderate deviations setting. In particular, this establishes the first claim of Proposition \ref{thm:moderate}.
	\begin{Prop}[Upper bound: right-deviation] \label{prop:upperboundx>1}
		For all $x\geq 1$,
		\begin{equation}\label{ineq:uppbound1}
        \limsup_{N\rightarrow\infty} \frac{1}{\log N}\log\P_N\left(\frac{x_\mrm{max}}{(\log N)^{1/ d}}\geq x\right) \leq -I_ d(x)\,.
		\end{equation}
        Furthermore, let $- d^{-1}<\gamma<1- d^{-1}$ and $x>0$. If we assume in addition $\varphi$ bounded, the following estimate holds:
			\begin{equation}\label{ineq:uppbound2}\limsup_{N\rightarrow\infty}\dfrac{1}{(\log N)^{1-\gamma- d^{-1}}}\log\mbb{P}_N\left(x_{\mrm{max}}\geq (\log N)^{1/ d}+\dfrac{x}{(\log N)^\gamma}\right)\leq -xd
            \end{equation}
	\end{Prop}

   
	\begin{Pro}
		By \cite[Proposition 2.1]{lambert2021poisson} (see also \cite[Proposition 3.6]{NakanoTrinh20} for the log case), there exists $C>0$, such that, setting $\rho_N$ the density of the first marginal given by
        \begin{equation}
            \label{ineq:wegner}
            \rho_N(u)\defi\int_{\R^{N-1}}p_N(u,x_2,\ldots,x_N)\diff x_2\ldots \diff x_N \leq C e^{-V(u)+2\alpha\vartheta(u)}\,.
        \end{equation}
        for all $u\in\R$, and $\vartheta \defi \log(1+|\cdot|)\mathbf{1}_{g=-\log}$,
        
		Thus, for all $N\geq 1$ and $x\geq 1$ (possibly depending on $N$), by a union bound and exchangeability we have:
		\begin{align*}
			\P_N\left(\frac{x_\mrm{max}}{(\log N)^{1/ d}}\geq x\right)&\leq N\P_N\left(\frac{x_1}{(\log N)^{1/ d}} \geq x\right) \\
			&\leq NC\int_{x(\log N)^{1/ d}}^{+\infty} \exp\left(-|u|^ d-\varphi(u)+2\alpha\vartheta(u)\right)\diff u \\
			&\leq N(\log N)^{1/ d}C\int_{x}^{+\infty} \exp\left(-v^ d\log N - \varphi\big((\log N)^{1/ d}v\big) + 2\alpha\vartheta\big((\log N)^{1/ d}v\big)\right)\diff v.
		\end{align*}
		Setting $\phi(x)=-\varphi(x)+2\alpha\vartheta(x)$ and since $\partial_v\left( v\mapsto e^{-(\log N)v^ d + \phi((\log N)^{1/ d}v)} \right)$ is negative for all $v\geq 1$ for $N$ large enough by Assumptions \ref{assumptions}, for all $t\geq x$
        \begin{align*}
            e^{-(\log N)t^ d + \phi((\log N)^{1/ d}t)} &= \frac{-1}{ d t^{ d-1}\log N-(\log N)^{1/ d}\phi'(t(\log N)^{1/ d})}\partial_t\left( e^{-(\log N)t^ d + \phi((\log N)^{1/ d}t)} \right)\\
            &\leq \dfrac{-2}{dt^{ d-1}\log N}\partial_t\left( e^{-(\log N)t^ d + \phi((\log N)^{1/ d}t)} \right).
        \end{align*}
        Thus, using that $ d\geq 1$, we obtain the following bound, valid for all $N$ large enough and $x\geq 1$:
        \begin{equation}
        \label{eq:bornesupprecise}
            \P_N\left( \frac{x_\mathrm{max}}{(\log N)^{1/ d}}\geq x \right)\leq  C'\frac{N(\log N)^{1/ d-1}}{x^{ d-1}}\exp\Big(-(\log N)x^ d + \phi((\log N)^{1/ d}x)\Big).
        \end{equation}
       Taking the logarithm and letting $N\rightarrow\infty$, yields \eqref{ineq:uppbound1} for fixed $x\geq 1$ by Assumption \ref{assumptions}. To establish \eqref{ineq:uppbound2} under the additional assumption that $\varphi$ is bounded, fix $x>0$ and use \eqref{eq:bornesupprecise} for $x_N=1+x/(\log N)^{1/ d + \gamma}$.
       Taylor expanding in \eqref{eq:bornesupprecise}, and using that by assumption, $\phi((\log N)^{1/ d}x_N)=O(\log\log N)$, we find
       \begin{align*}
           \P_N\left( \frac{x_\mathrm{max}}{(\log N)^{1/ d}}\geq x \right) &\leq \frac{C'(\log N)^{1/ d-1}}{ x_N^{ d-1}} N\exp\Big(-(\log N)\big(1+ xd/(\log N)^{1/ d+\gamma}\big)+o((\log N)^{1-1/ d - \gamma})\Big)\\
           &\leq \exp(- xd (\log N)^{1-1/ d - \gamma}+o((\log N)^{1-1/ d - \gamma})),
       \end{align*}
       yielding \eqref{ineq:uppbound2}.
	\end{Pro}

    \subsubsection{Exponential tightness}
	\begin{Cor}[Exponential tightness]\label{cor:exptight}
		The law of $(\overline{x_{\max}})_N$ is exponentially tight, \textit{i.e.} there exists a sequence $(K_M)_{M\geq 1}$ of compacts of $\R$ such that:
		\begin{equation}
			\label{eq:expTight}
		\lim_{M\rightarrow\infty}\limsup_{N\rightarrow\infty} \frac{1}{\log N}\log \P_N\left( \frac{x_{\max}}{(\log N)^{1/ d}} \in K_M^\mathsf{c}  \right) = -\infty\,.
		\end{equation}
	\end{Cor}
	
	\begin{Pro}
		Indeed let $M\geq1$, 
		\begin{equation*}
	\limsup_{N\rightarrow\infty}\frac{1}{\log N}\log \P_N\Big( \frac{x_{\max}}{(\log N)^{1/ d}}\in [-M,M]^\mathsf{c}  \Big) = \max\{a_M,b_M\}
		\end{equation*}
		where
		$$ a_M \defi  \limsup_{N\rightarrow\infty}\frac{1}{\log N}\log\P_N\Big( \frac{x_{\max}}{(\log N)^{1/ d}}<-M \Big) $$
		and
		$$ b_M \defi  \limsup_{N\rightarrow\infty}\frac{1}{\log N}\log\P_N\Big( \frac{x_{\max}}{(\log N)^{1/ d}}>M \Big)\,.$$

        By Proposition \ref{prop:upperboundx>1}, and by definition of $I_d$, we have $$b_M\leq-I_ d(M)\underset{M\to +\infty}{\longrightarrow}-\infty.$$ 
        Similarly, adapting the proof of Proposition \ref{prop:upperboundx>1} for $x_{\min}$, one shows under the same assumptions that the following bound for $x\leq -1$ holds:
        $$ \limsup_{N\rightarrow\infty} \frac{1}{\log N}\log \P\Big(\frac{x_{\min}}{(\log N)^{1/d}}\leq x\Big)\leq -I_d(x).$$ 
        Using that $x_{\min}\leq x_{\max}$, we conclude that $\lim_{M\to \infty}a_M = -\infty$,
        proving \eqref{eq:expTight} with $K_M=[-M,M]$.
	\end{Pro}
	
\subsubsection{Upper bound: left deviation.}\label{subsec:left}
We now prove an upper bound on the probability of a left-deviation of $x_{\max}$. Our strategy relies on a comparison with the edge of a system of iid particles distributed according to $\bigotimes^N\mu_{\mrm{eq}}$, where $\mu_\mrm{eq}$ is the unique minimizer of \eqref{def:energie}. In the latter case, deriving the precise asymptotic of a left deviation follows from elementary computations.
Prior to that, we establish the following controls on $U^{\rho_\mrm{eq}}$, defined in \eqref{eqmeasure:exp}.

\begin{Lem}\label{lem:U rho eq ' bornée}
 The function $U^{\rho_\mrm{eq}}$ is differentiable on $\R$ and its derivative is a bounded function on $\R$. Furthermore, as $|x|\rightarrow\infty$
 $$U^{\rho_\mrm{eq}}(x)=\log(1+|x|)\mbf{1}_{g=-\mrm{log}}+\mbf{1}_{g=-\mrm{log}}O(x^{-1})+\mbf{1}_{g=|\cdot|^{-s}}O(x^{-s}).$$
\end{Lem}

\begin{Pro}
   The proof is done for the log case in \cite[Lemma 2.3]{DwoMemin} so we only consider the Riesz case here. From \eqref{eqmeasure:charac} and the fact that $U^{\rho_\mrm{eq}}(x)\geq0$, we see that for all $x$, $\rho_\mrm{eq}(x)\leq C e^{-V(x)}$ for some $C>0$ which ensures that $\rho_\mrm{eq}\in L^2(\R)$. Furthermore, in the sense of distributions:
     $$\partial_xU^{\rho_\mrm{eq}}=\partial_x\left(\rho_\mrm{eq}\ast|\cdot|^{-s}\right)=-s\rho_\mrm{eq}\ast\mathrm{p.v}\left(\dfrac{\mrm{sgn}(\cdot)}{|\cdot|^{s+1}}\right).$$
     where $\mrm{p.v}$ denotes the Cauchy principal value distribution. The distribution in the RHS defines a function in $L^1_\mrm{loc}(\R)$, thus $U^{\rho_\mrm{eq}}$ is a differentiable function and
   $$[U^{\rho_\mrm{eq}}]'(x)=-s\mcal{H}_s[\rho_\mrm{eq}](x)\defi-s\fint_\R \dfrac{\rho_\mrm{eq}(y)\mrm{sgn}(x-y)}{|x-y|^{s+1}}\diff y.$$
   Furthermore, for all $x\in\R$:
   $$\mcal{H}_s[\rho_\mrm{eq}](x)=\bigg(\fint_{|y|<|x|/2}+\fint_{|y|>|x|/2}\bigg)\dfrac{\rho_\mrm{eq}(y)\mrm{sgn}(x-y)}{|x-y|^{s+1}}\diff y.$$
   On $\{|y|<|x|/2\}$, $|x-y|>|x|/2$, while on $\{|y|>|x|/2\}$, because of Assumptions \ref{assumptions}, there exists $c>0$ such that $\rho_\mrm{eq}(y)\leq e^{-c|x|^d}$ for $|x|$ large enough. This yields:
   $$|\mcal{H}_s[\rho_\mrm{eq}](x)|\leq\dfrac{2}{|x|^{s+1}}+2e^{-c|x|^d}\int_{|y|>|x|/2}\dfrac{\diff y}{|y|^{s+1}}=O(x^{-s-1}).$$
   Thus after integrating, we obtain that $U^{\rho_\mrm{eq}}(x)=O(x^{-s})$ at infinity. It only remains to show that $\mcal{H}_s[\rho_\mrm{eq}]$ is bounded on all compact sets. For that, let $r\in\R$, we introduce the $r$-th Sobolev space:
   $$H_r(\R)=\left\{u\in\mcal{S}'(\R),\|u\|_{H^{r}(\R)}^2\defi\int_{\R}(1+|\xi|)^{2r}\cdot|\mcal{F}[u](\xi)|^2\diff\xi<+\infty\right\},$$
   where $\mcal{S}'(\R)$ denotes the space of Schwartz distributions and where $\mcal{F}[u]\defi\int_\R u(x)e^{\ii x\cdot}\diff x$. Taking the Fourier transform, we obtain that there exists $c_s\in\R$:
   $$\mcal{F}\left[\mcal{H}_s[\rho_\mrm{eq}]\right](\xi)=\ii c_s\mrm{sgn}(\xi)|\xi|^{s}\mcal{F}\left[\rho_\mrm{eq}\right](\xi),$$
   thus
   $$\|\mcal{H}_{s}[\rho_\mrm{eq}]\|_{H^{-s}(\R)}^2\leq c_s^2\int_{\R}|\mcal{F}[\rho_\mrm{eq}](\xi)|^2\diff\xi\leq c_s^2\|\rho_\mrm{eq}\|_{L^2(\R)}^2.$$
   We thus, conclude that $\mcal{H}_{s}[\rho_\mrm{eq}]\in H^{-s}(\R)$ and thus $\rho_\mrm{eq}'=(-V'-s\mcal{H}_s[\rho_\mrm{eq}])\rho_\mrm{eq}\in H^{-s}(\R)$. We conclude from that and $\rho_\mrm{eq}\in L^2(\R)=H^0(\R)$ that $\rho_\mrm{eq}\in H^{1-s}(\R)$. Thus, repeating the same argument, we have:
   $$\|\mcal{H}_{s}[\rho_\mrm{eq}]\|_{H^{1-2s}(\R)}^2\leq c_s^2\int_{\R}\cdot(1+|\xi|)^{2(1-s)}\cdot|\mcal{F}[\rho_\mrm{eq}](\xi)|^2\diff\xi\leq c_s^2\|\rho_\mrm{eq}\|_{H^{1-s}(\R)}^2.$$
   We deduce again that, $\rho_\mrm{eq}'\in H^{0\wedge(1-2s)}(\R)$ and thus that $\rho_\mrm{eq}\in H^{1\wedge(2-2s)}(\R)$. Iterating the argument yields that $\rho_\mrm{eq}\in H^{1\wedge (n-ns)}(\R)= H^1(\R)\subset\mcal{C}^0(\R)$ for $n>1/(1-s)$. This concludes the proof.$\hfill\square$\end{Pro}

We will use the following distance on $\mathcal{M}_1(\R)$.
\begin{Def}
    Define the distance (possibly infinite) on $\mathcal{M}_1(\R)$ by
\begin{align*}
D(\mu,\nu)&\defi\Big(\iint_{\R^2}g(x-y)\diff(\mu-\nu)(x)\diff(\mu-\nu)\Big)^{1/2}\\
&= \Big(\int_{\R}|t|^{\mathsf{\mathbf{s}}-1}|\mathcal{F}(\mu-\nu)(t)|^2\diff t\Big)^{1/2},
\end{align*}
where we denoted $\mathsf{\mathbf{s}}= s$ if $g=|x|^{-s}$, and $0$ if $g=-\log$, and wrote the Fourier transform of a signed measure $\sigma$ by $\mathcal{F}[\sigma](x)=\int_\R e^{\mathrm{i}tx}\diff\sigma(t).$
\end{Def}

We then have the following identity, which comes as a consequence of \eqref{eqmeasure:charac} and is a direct extension of \cite[Lemma 3.3]{DwoMemin}.
\begin{Lem}\label{lem:distance}
For any $\mu\in \mathcal{M}_1(\R)$ with finite energy $\mathcal{E}(\mu)$, we have
$$ \mathcal{E}(\mu)-\mathcal{E}(\mu_{\mathrm{eq}}) = \alpha D^2(\mu,\mu_{\mathrm{eq}})+\int_{\R}\log\frac{\diff \mu}{\diff \mu_{\mathrm{eq}}}(x)\diff\mu(x).$$
\end{Lem}

We are now able to compare the law of the edge of $\mbb{P}_N$ to the one of $\mu_\mrm{eq}^{\otimes N}$.

\begin{Lem}\label{lem:iidcomparison} Let $x\in(0,1)$, then
    \begin{equation*}
         \P_N\left( \overline{x_{\max}}<(1-x)^{1/ d}\right) \leq e^{\mfrak{R}_N}\mu_\mrm{eq}^{\otimes N}\left( \overline{x_{\max}}<(1-x)^{1/ d} \right)
    \end{equation*}
    where $\mfrak{R}_N=O(\log N)$ for the case $g=-\log$ and $\mfrak{R}_N=O(N^{\frac{2s}{1+s}})$ for the case $g=|.|^{-s}$.
\end{Lem}

\begin{Pro}We prove this bound by bounding the density $p_N$ pointwise. To do so, we use energy estimates. Since $\mcal{E}(\mu_N)=+\infty$ for all $\underline{x}\in\R^N$, we need to approximate $\mu_N$ by an absolutely continuous measure with respect to Lebesgue measure. As in \cite[Section 3.1]{DwoMemin}, we introduce $\widetilde{\mu}_N$, the regularization of $\mu_N$ by setting, for any $\underline{x} \in \R^N$, the following: we first apply a permutation $\sigma$ which orders the particles $x_{\sigma(1)}\leq\hdots\leq x_{\sigma(N)}$ and then we define $\underline{y}\in\R^N$ as follows:
    $$y_1 \defi x_{\sigma(1)},\hspace{1cm}\forall k\in\llbracket0,N-1\rrbracket,\,\,\,y_{k+1}\defi y_k + \max\{x_{\sigma(k+1)}-x_{\sigma(k)}, N^{-b}\}.$$
    Finally let $a>0$, we define $\lambda_{N^{-a}}$ to be the uniform measure on $[0,N^{-a}/2]$ and $\widetilde{\mu}_N$ as the convolution between $\lambda_{N^{-a}}$ and $\mu_N$, namely:
    \begin{equation}
    \label{eq:tildeMu}
    \widetilde{\mu}_N\defi\lambda_{N^{-a}} \ast \frac{1}{N}\sum_{i=1}^N \delta_{y_{\sigma^{-1}(i)}}=\lambda_{N^{-a}} \ast\frac{1}{N}\sum_{i=1}^N \delta_{y_{i}}.
\end{equation} When $g=-\log$, we take $a=5$ and $b=3$, while for the case $g=|.|^{-s}$, we take $a=b=\frac{2}{1+s}$.

As in the proof of \cite[Proposition 1.6]{DwoMemin}, it holds that $\mcal{Z}_N\geq e^{-N\mcal{E}(\mu_{\mrm{eq}})+c}$ for some constant $c$ depending on $V,P$ and $s$ in the Riesz case. Using that for all $i\neq j$ one has $|x_i-x_j|\leq |y_i-y_j|$ and that $g$ is non-increasing on $\R_+$, we have:
$$\log p_N(\underline{x})\leq N\mcal{E}(\mu_{\mrm{eq}})-\dfrac{\alpha}{N}\sum_{i\neq j}g(y_i-y_j)-\sum_{i=1}^NV(x_i)+O(1).$$
For all $x\neq 0$ and $|h|\leq\dfrac{|x|}{2}$ we have $|g(x+h)-g(x)|\leq C|g'(x)|.|h|$ for a uniform constant $C>0$. Thus, we have: for all $u,v\in[0,N^{-a}/2]$, and $i\neq j$, 
$$|g(y_i+u-y_j-v)-g(y_i-y_j)|\leq C|u-v|.|g'(y_i-y_j)|\leq C.N^{-a}|g'(N^{-b})|\defi A_N.$$
Above, we used that $|y_i-y_j| \geq N^{-b}$  for all $i\neq j$. $A_N$ is a $O(N^{-2})$ in the log case and $O(N^{\frac{2s}{1+s}})$ in the Riesz case. Now summing over $i\neq j$ and integrating against $\diff u\diff v$, we obtain: 
\begin{multline}
    -\dfrac{\alpha}{N}\sum_{i\neq j}g(y_i-y_j)\leq 2\alpha N A_N-\dfrac{\alpha}{N}\sum_{i\neq j}\iint_{[0,N^{-a}/2]^2}g(y_i-y_j+u-v)\diff \lambda_{N^{-a}}(u)\diff\lambda_{N^{-a}}(v)
    \\\leq 2\alpha N A_N-\alpha N\iint_{\R^2}g(y-z)\diff\widetilde{\mu}_N(y)\diff\widetilde{\mu}_N(z)+B_N
\end{multline}
where the diagonal term $$B_N\defi \alpha\iint_{[0,N^{-a}/2]^{2}}g(u-v)\diff\lambda_{N^{-a}}(u)\diff\lambda_{N^{-a}}(v)$$
is a $O(\log N)$ in the log case and $O(N^{\frac{2s}{1+s}})$ in the Riesz case, in particular $NA_N=O(B_N)$. Recalling the definition of $\mcal{E}$ \eqref{def:energie}, we obtain thanks to Lemma \ref{lem:distance}:
\begin{align*}
    \log p_N(\underline{x})&\leq- N\big(\mcal{E}(\widetilde{\mu}_N)-\mcal{E}(\mu_\mrm{eq})\big)++N\int_\R\big[V(y)+\log\dfrac{\diff\widetilde{\mu}_N}{\diff y}(y)]\diff\widetilde{\mu}_N(y)-\sum_{i=1}^NV(x_i)+O(B_N)
    \\&\leq-N\alpha D^2(\widetilde{\mu}_N,\mu_{\mathrm{eq}})
    +N\int_\R\big[V(y)+\log\rho_\mathrm{eq}(y)]\diff\widetilde{\mu}_N(y)-\sum_{i=1}^NV(x_i)+O(B_N).
\end{align*}
Finally we use \eqref{eqmeasure:charac} and the fact that $D^2(\widetilde\mu_N,\mu_\mathrm{eq})\geq 0$
to deduce that
\begin{equation*}
    \log p_N(\underline{x})\leq -2\alpha N\int_{\R}U^{\rho_\mrm{eq}}(y)\diff(\widetilde{\mu}_N-\mu_N)(y)+\sum_{i=1}^N\log\rho_\mrm{eq}(x_i)+O(B_N).
\end{equation*}
Finally, we use Lemma \ref{lem:U rho eq ' bornée} to deduce that
\begin{multline}
    \Big|N\int_{\R}U^{\rho_\mrm{eq}}(x)\diff(\widetilde{\mu}_N-\mu_N)(x)\Big|\leq\Big\|\dfrac{\diff}{\diff x}U^{\rho_\mrm{eq}}\Big\|_\infty\Big(\sum_{i=1}^N|y_i-x_i|+N\int_0^{N^{-a}}u\diff\lambda_{N^{-a}}(u)\Big)
    \\\leq\Big\|\dfrac{\diff}{\diff x}U^{\rho_\mrm{eq}}\Big\|_\infty O\big([N^{\frac{2s}{1+s}}+N^{\frac{-1+s}{1+s}}]\mbf{1}_{g=|.|^{-s}}+[N^{-1}+N^{-4}]\mbf{1}_{g=-\log}\big)
\end{multline}
where we used that $|x_i-y_i|\leq(i-1)N^{-b}$.
These estimations lead to the following bound:
$$p_N(\underline{x})\leq e^{O(B_N)}.\prod_{i=1}^N\rho_\mrm{eq}(x_i).$$ This concludes the proof.
\end{Pro}

The following proposition deduces from Lemma \ref{lem:iidcomparison} an upper bound on left deviations.

\begin{Prop}[Upper bound: left-deviation]\label{prop:LUP}
Under Assumptions \ref{assumptions}, assume that $\varphi$ is bounded. Set $x_0 = 0$ for $g=-\log$, and $x_0=\frac{2s}{1+s}$ for $g=|.|^{-s}$.
    For all $x_0<x<1$,
    \begin{equation}\label{ineq:Luppbound1}
        \limsup_{N\rightarrow\infty}\dfrac{1}{\log N}\log\log \Big(\P\Big(x_{\max}\leq((1-x)\log N)^{1/ d}\Big)^{-1}\Big)\geq x.
    \end{equation}
    Furthermore, for $g=-\log$ and for all $\gamma$ such that $- d^{-1}<\gamma<1- d^{-1}$, we have for all $x>0$:
           \begin{equation}\label{ineq:Luppbound2}
			\limsup_{N\rightarrow\infty}\dfrac{1}{(\log N)^{1-\gamma- d^{-1}}}\log\log \Big(\mbb{P}_N\Big(x_{\mrm{max}}\leq (\log N)^{1/ d}-\dfrac{x}{(\log N)^\gamma}\Big)^{-1}\Big)\geq xd
		\end{equation}
\end{Prop}
\begin{Rem}
In particular, \eqref{ineq:Luppbound1} establishes Proposition \ref{prop:preciseupper}, and \eqref{ineq:Luppbound2} establishes the second claim of Proposition \ref{thm:moderate}. Together with Proposition \ref{prop:upperboundx>1}, this completes the proof of Proposition \ref{thm:moderate}.
\end{Rem}

\begin{Pro}
We first prove that there exists $C>0$ such that as $x\rightarrow+\infty$:
$$\mu_\mrm{eq}\big((x,+\infty)\big)=C\dfrac{(1+|x|)^{2\alpha \mbf{1}_{g=-\log}}}{V'(x)}e^{-V(x)}\big(1+o(1)\big).$$
Let $x>0$, an integration by parts yields:
$$\mu_\mrm{eq}\big((x,+\infty)\big)=\int_x^{+\infty}\dfrac{\rho_\mrm{eq}}{\rho_\mrm{eq}'}(t)\cdot\rho_\mrm{eq}'(t)\diff t=\dfrac{-\rho_\mrm{eq}}{\rho_\mrm{eq}'}(x)\cdot\rho_\mrm{eq}(x)-\int_x^{+\infty}\left(\dfrac{\rho_\mrm{eq}}{\rho_\mrm{eq}'}\right)'(t)\cdot\rho_\mrm{eq}(t)\diff t.$$
Since, by Assumption \ref{assumptions}, $|V'(x)|\rightarrow+\infty$ and $V''/V'\to 0$ at $\infty$, one checks, as in \cite[Lemma 2.2]{DwoMemin}, using \eqref{eqmeasure:charac} and Lemma \ref{lem:U rho eq ' bornée} that $\rho_\mathrm{eq}$ is indeed $\mathcal{C}^2$ outside a compact set, and that
$\left(\dfrac{\rho_\mrm{eq}}{\rho_\mrm{eq}'}\right)'(x)=o(1)$ as $x\rightarrow\infty$, and thus, by comparing it with the first integral above, the second integral in the RHS can be seen to be a $o\Big(\mu_\mrm{eq}\big((x,+\infty)\big)\Big)$. We thus obtain, by \eqref{eqmeasure:exp}:
\begin{align*}
    \mu_\mrm{eq}\big((x,+\infty)\big)&=\dfrac{-\rho_\mrm{eq}}{\rho_\mrm{eq}'}(x)\cdot\rho_\mrm{eq}(x)\big(1+o(1)\big)
\\&=\dfrac{1}{V'(x)\big(1+o(1)\big)}\exp\left(-V(x)+2\alpha \log(1+|x|)\mbf{1}_{g=-\log}+\lambda_\mrm{eq}+o(1)\right).
\end{align*}
This concludes the first step.
For \eqref{ineq:Luppbound1}, by direct computation, as $N$ goes to infinity, we have, denoting $u_N=((1-x)\log N)^{1/ d}$:
\begin{align*}
    \mu_\mrm{eq}^{\otimes N}\left(x_{\max}\leq u_N\right)&=\left[1-\mu_\mrm{eq}\left(u_N,+\infty)\right)\right]^N
    \\&=\exp\left(N\log\Big[1-C\dfrac{(1+|u_N|)^{2\alpha \mbf{1}_{g=-\log}}}{V'(u_N)}e^{-V(u_N)}\big(1+o(1)\big)\Big]\right)
    \\&\leq \exp\left(-N\exp[-V(u_N)+o(V(u_N)]\right)
    \\&\leq \exp\left(-N\exp[-(1-x)\log N+o(\log N)]\right).
\end{align*}
Thus, using Lemma \ref{lem:iidcomparison}, for the log-case, we obtain:
\begin{equation}\label{eq:prop2.5log}
    \log\mbb{P}_N\left(x_{\max}\leq((1-x)\log N)^{1/ d}\right)\leq-N^{x+o(1)}+O(\log N)
\end{equation}
which allows us to conclude. Similarly, for the Riesz-case, we obtain:
\begin{equation}\label{eq:prop2.5riesz}
   \log\mbb{P}_N\left(x_{\max}\leq((1-x)\log N)^{1/ d}\right)\leq-N^{x+o(1)}+O(N^{\frac{2s}{1+s}}).
\end{equation}
Thus for any $x>\frac{2s}{1+s}$, we can deduce the desired bound. To obtain \eqref{ineq:Luppbound2}, we use the exact same arguments with $u_N=(\log N)^{1/ d}-\dfrac{x}{(\log N)^\gamma}$. Taylor expanding,
\begin{multline*}
    \log N-V(u_N)=\log N-\Big((\log N)^{1/ d}-\dfrac{x}{(\log N)^\gamma}\Big)^ d+\varphi(u_N) \\
    =d\dfrac{  x}{(\log N)^{\gamma-1-1/ d}}+o\Big(\dfrac{ x}{(\log N)^{\gamma-1-1/ d}}\Big)+\varphi(u_N).
\end{multline*}
We thus obtain:
$$\dfrac{1}{(\log N)^{1+1/ d-\gamma}}\log\log\mbb{P}_N\left(x_{\max}\leq((1-x)\log N)^{1/ d}\right)^{-1}\geq xd+\dfrac{\varphi\left(u_N\right)}{(\log N)^{1+1/ d-\gamma}}+O(\log N).$$
We then conclude by boundedness of $\varphi$.
\end{Pro}

In particular, we have the immediate following corollary which gives the upper bound on left-deviations.

\begin{Cor}\label{cor:left_upperbound} For all $x_0<x<1$,
    $$\limsup_{N\rightarrow\infty}\dfrac{1}{\log N}\log\P\left(\dfrac{x_{\max}}{(\log N)^{1/ d}}\leq(1-x)^{1/ d}\right)=-\infty=-I_d(x).$$
\end{Cor}

\begin{Pro}
    It is a direct consequence of \eqref{eq:prop2.5log} and \eqref{eq:prop2.5riesz}.
\end{Pro}

    \subsection{Lower bound}
For the following, let $\upbeta_N>0$ (not necessarily $2\alpha/N$), we set $\P_N^{\upbeta_N}$ to be the probability measure on $\R^N$ with pdf with respect to the Lebesgue measure given by: 
$$\dfrac{1}{\mcal{Z}_N^{\upbeta_N}}\prod_{i<j}^{N}e^{-\upbeta_Ng(x_i-x_j)}\cdot\prod_{i=1}^{N}e^{-V(x_i)}.$$

Under Assumptions \ref{assumptions}, if $N\upbeta_N\to 2\alpha$, as shown in \cite[Theorems 3.2, 3.4]{lambert2021poisson} the maximum of the configuration $\{x_i\}$ under $\P_N^{\upbeta_N}$ fluctuates around $E_N \sim (\log N)^{1/ d}$, both in the Riesz and the $\log$ cases.
The edge $E_N$ is identified in the Riesz case through:
\begin{equation}
    \label{eq:edgeRiesz}
    \frac{Ne^{-V(E_N)}}{e^{\lambda_\mathrm{eq}}V'(E_N)} \underset{N\to \infty}{\longrightarrow} 1,
\end{equation}
where $\lambda_\mathrm{eq}$ is the constant of \eqref{eqmeasure:charac}. In the $\log$-gas case, the equation for $E_N$ is similar and given by:
\begin{equation}
    \label{eq:edgeLog}
    \frac{Ne^{-V(E_N)+N\upbeta_N \log E_N}}{e^{\lambda_\mathrm{eq}}V'(E_N)} \underset{N\to \infty}{\longrightarrow} 1.
\end{equation}

\begin{Lem}[PPP at the edge]
\label{lem:lambertedge}
    Assume that $N\upbeta_N\tend{N\rightarrow\infty}2\alpha$ and  $(E_N)_N$ satisfies respectively \eqref{eq:edgeRiesz} in the Riesz case ($g=|\cdot|^{-s}$), and \eqref{eq:edgeLog} in the log gas case $(g=-\log)$. Then,
    \begin{equation}\label{lem:PPPedge}
        E_N \underset{N\to \infty}{\sim} (\log N)^{1/ d},
    \end{equation}and under $\P_N^{\upbeta_N}$ the sequence of random measures $(\Xi_N)_N$ defined by
    $$ \Xi_N \defi \sum_{j=1}^N \delta_{V'(E_N)(x_j-E_N)} $$
    converges in distribution towards an inhomogeneous Poisson point process with intensity $x\mapsto e^{-x}$, meaning that for all Borel set $A\subset \R$ such that $\int_A e^{-x}\diff x<\infty$, the random variable $\Xi_N(A)$ converges towards a Poisson random variable with mean $\int_A e^{-x}\diff x$.
\end{Lem}

\begin{Pro}
The only point that needs to be checked is the asymptotic \eqref{lem:PPPedge} of $E_N$. It is established for the $\log$ gas in \cite[Theorem 4.1]{DwoMemin}, and the proof for the Riesz gas is identical.
\end{Pro}

\begin{Lem}\label{lem:forlowerbound} Assume that $N\upbeta_N\tend{N\rightarrow\infty}2\alpha$, there exists $p\in(0,1)$ such that for $N$ large enough:
    \begin{equation}
    \label{eq:p>0}
    \P_N^{\upbeta_N}\left(\forall j,\ x_j<E_N\right) \geq p.
\end{equation}
Furthermore, there exists $C>0$ such that:
$\dfrac{\mcal{Z}_{N-1}^{\upbeta_N}}{\mcal{Z}_{N}^{\upbeta_N}}\tend{N\rightarrow\infty}C$. Finally, the random variable $\overline{x_{\max}}$ converges in probability to $1$.
\end{Lem}

\begin{Pro} As a consequence of Lemma \ref{lem:lambertedge}:
$$ \P_N^{\upbeta_N}\left(\forall j,\ x_j<E_N\right) = \P_N^{\upbeta_N}\big(\Xi_N(0,\infty\big) =0) \tend{N\rightarrow\infty}\exp\Big(-\int_{(0,\infty)}e^{-t}\diff t\Big)= e^{-1}.$$
This establishes the first point.
    The second point comes from \cite[Corollary 2.4]{lambert2021poisson}. The last point is a consequence of \cite[Theorem 3.2, 3.4]{lambert2021poisson}, showing that $dE_N^d(\overline{x_{\max}}-1)$ converges in law towards a Gumbel random variable.
\end{Pro}

As a consequence, we have:

\begin{Prop}[Lower bound]
		\label{pro:lowerbound}
		Let $x\in\R$, the following bound holds:
		\begin{equation*}
			\lim_{\delta\downarrow0} \liminf_{N\rightarrow\infty}\dfrac{1}{\log N}\log \mbb{P}_N\Big(\overline{x_{\max}}\in(x-\delta,x+\delta)\Big)\geq-I_ d(x).
		\end{equation*}
	\end{Prop}
	
	\begin{Pro}When $x<1$, $I_ d(x)=+\infty$ and there is nothing to prove, while for $x=1$, by Lemma \ref{lem:forlowerbound}, we have
    $$\dfrac{1}{\log N}\log \mbb{P}_N\left(\Big|\overline{x_{\max}}-1\Big|<\delta\right)\tend{N\rightarrow\infty}0$$
    and we conclude by noticing that $-I_ d(1)=0$. Let $x>1$ and $\delta>0$ small enough. To show that the above inequality is true, we restrict ourselves to the configurations such that there is exactly one particle abnormally far away from the finite $N$ expected support (namely $[-E_{N-1},E_{N-1}]$) and the $N-1$ remaining particles are located in this expected support. Let
		$$ I_N \defi  (-\infty,E_{N-1}].$$ 
		By noticing that 
		$$\biguplus_{i=1}^N\Big\{\Big|\dfrac{x_i}{(\log N)^{1/ d}}-x\Big|<\delta,\;\forall j\neq i,\ x_j\in I_N \Big\} \subset\Big\{\Big|\dfrac{x_{\mrm{max}}}{(\log N)^{1/ d}}-x\Big|<\delta\Big\},$$
		and setting $\upbeta_N=\beta_N=2\alpha/N$, we have by exchangeability, 
		\begin{multline*}
			\mbb{P}_N\Big(\Big|\dfrac{x_{\mrm{max}}}{(\log N)^{1/ d}}-x\Big|<\delta\Big)\geq N\mbb{P}_N\Big(\Big|\dfrac{x_1}{(\log N)^{1/ d}}-x\Big|<\delta,\;\forall j>1,\ x_j\in I_N\Big)
			\\=N\frac{\mcal{Z}_{N-1}^{\upbeta_N}}{\mcal{Z}_N^{\upbeta_N}}\int_{(x-\delta)(\log N)^{1/ d}}^{(x+\delta)(\log N)^{1/ d}} e^{-V(x_1)}\int_{I_N^{N-1}}e^{-\upbeta_N\sum_{i=2}^Ng(x_1-x_i)}\diff\P^{\upbeta_N}_{N-1}(x_2,\ldots,x_N)\diff x_1
			\\= 2\delta(\log N)^{1/ d}N\dfrac{\mcal{Z}_{N-1}^{\upbeta_N}}{\mcal{Z}_{N}^{\upbeta_N}}\cdot\P_{N-1}^{\upbeta_N}\left(\forall j\geq 2, x_j \in I_N\right)\cdot\mcal{I}
		\end{multline*}
        where
        \begin{align*}
            \mcal{I}\defi&\int_{x-\delta}^{x+\delta} \exp\left(-V(u(\log N)^{1/ d})\right)
            \\&\quad\quad\times\int_{I_N^{N-1}}\exp\left(-\upbeta_N\sum_{i=2}^Ng\left(u(\log N)^{1/ d}-x_i\right)\right)\diff\mu(u,x_2\ldots,x_N)
        \end{align*}
		and where the probability measure $\mu$ is supported on $[x-\delta,x+\delta]\times I_N^{N-1}$ and defined by $$\diff\mu(u,\ldots,x_N)\defi \mbf{1}_{|u-x|<\delta}\dfrac{\diff u}{2\delta}\otimes\mbf{1}_{\cap_j x_j\in I_N}\frac{\diff\mbb{P}_{N-1}^{\upbeta_N}(x_2,\dots,x_N)}{\P_{N-1}^{\upbeta_N}\left(\forall j\geq 2, x_j \in I_N\right)}.$$ 
		By Lemma \ref{lem:forlowerbound}, taking the logarithm and applying Jensen's inequality, we obtain: 
		\begin{multline*}
			\log \mbb{P}_N\left(\Big|\dfrac{x_{\mrm{max}}}{(\log N)^{1/ d}}-x\Big|<\delta\right)\geq \log N+ o_N(\log N)
            \\ +
			\int_{\R^N}\left[-u^{ d}\log N-\varphi(u(\log N)^{1/ d}) -\upbeta_N\sum_{i=2}^Ng\left(u(\log N)^{1/ d}-x_i\right)\right]\diff\mu(u,x_2,\ldots,x_N)
			\\\geq \log N+ o_N(\log N)\\-\log N (x+\delta)^ d-\int_{x-\delta}^{x+\delta} \varphi(u(\log N)^{1/ d})\frac{\diff u}{2\delta}
			-\int_{\R^N}\upbeta_N\sum_{i=2}^Ng\left(u(\log N)^{1/ d}-x_i\right)\diff\mu(u,x_2,\ldots,x_N).
            \end{multline*}
Since by Assumptions \ref{assumptions} we have for $\delta>0$ fixed,
        $$ \frac{1}{2\delta}\int_{x-\delta}^{x+\delta} \varphi(u(\log N)^{1/ d})\diff u = o_N(\log N),$$
we thus have 
\begin{multline*}
			\dfrac{1}{\log N}\log \mbb{P}_N\left(\dfrac{x_\mrm{max}}{(\log N)^{1/ d}}\in(x-\delta,x+\delta)\right)
            \geq -I_ d(x+\delta)+o_N(1)\\
            -\dfrac{1}{\log N}\int_{[x-\delta,x+\delta]\times I_N^{N-1}}\upbeta_N\sum_{i=2}^Ng\left(u(\log N)^{1/ d}-x_i\right)\diff\mu(u,x_2,\ldots,x_N).
\end{multline*}
For $\delta$ small enough, it holds $x-1-2\delta >0$, and for $N$ big enough $\frac{E_{N-1}}{(\log N)^{1/ d}} < 1+\delta$. Hence for all $u\in (x-\delta,x+\delta)$, $y\in I_N$ we have $|u(\log N)^{1/ d}-y| \geq (x-1-2\delta)(\log N)^{1/ d}$, therefore since $-g$ is increasing:
$$ -\hspace{-0.3cm}\int\limits_{[x,x+\delta]\times I_N^{N-1}}\upbeta_N\sum_{i=2}^Ng\left(u(\log N)^{\frac1d}-x_i\right)\diff\mu(u,x_2,\ldots,x_N) \geq -(N-1)\upbeta_Ng\left((x-1-2\delta)(\log N)^{\frac1d}\right).$$
The RHS being $o_N(\log N)$ in both Riesz and log cases, we conclude that:
	$$\liminf_{N\rightarrow\infty}\dfrac{1}{\log N}\log \mbb{P}_N\left(\dfrac{x_\mrm{max}}{(\log N)^{1/ d}}\in(x-\delta,x+\delta)\right)\geq 1-(x+\delta)^ d.$$
    Letting $\delta\rightarrow0$ yields the conclusion.
\end{Pro}

\section{Matrix interpretation: proof of Theorem \ref{thm:LDPtridiag_intro}}\label{sec:matrix}
    In this section, we consider random matrix models of symmetric tridiagonal matrices with independent entries having Gaussian tails, and establish Theorem \ref{thm:LDPtridiag_intro}.

The following lemma will allow us to deduce the LD lower bound, and to control the spectrum of a matrix with small entries. In the following, all the matrices are of size $N\times N$, and we take the indices of the entries modulo $N$. In the nonperiodic model $T$, we consider $b_N=0$.
    \begin{Lem}\label{lem:inclusion_bourrine}Let $J=T$ or $J=T_\mathsf{per}$, for arbitrary $(a_i)_{1\leq i \leq N},(b_i)_{1\leq i \leq N}$. Let 
    $$ \rho(J)=\max_{\lambda\in \mathrm{Sp}(J)}|\lambda|$$
    be the spectral radius of $J$. Then,
    $$ \max_{1\leq i\leq N}(a_i) \leq \lambda_\mathrm{max}(J) \leq \rho(J) \leq \max_{1\leq i\leq N} \left(|a_i|+|b_{i-1}|+|b_i|\right).$$
    \end{Lem}
    
		
		\begin{Pro}
        For the lower bound, write $J=PDP^\top$ with $P$ orthogonal and $D=\mathrm{diag}(\lambda_1,\ldots,\lambda_N)$. Then, for any $1\leq i \leq N$, it holds that $a_i = J_{ii} = \displaystyle\sum_{j=1}^N p_{i,j}^2 \lambda_j \leq  \lambda_\mathrm{max}$.\\
	For the upper bound on $\rho(J)$, recall that the spectrum of any $N\times N$ matrix is contained in the union of the Gerschgorin disks: 
			$$\mathrm{Sp}(J)\subset \bigcup_{i=1}^N \overline{B}\left(J_{ii},\sum_{j\neq i} |J_{ij}|\right)\,, $$
			where $\overline{B}(x,\delta)=\{y\in \R:|x-y|\leq \delta\}$. Using the tridiagonal structure of $J$, we deduce that for some $1\leq i\leq N$,
            \begin{equation}
                \label{eq:bornelambda}
                |\lambda_{\max}| \leq |\lambda_{\max} - a_i| + |a_i| \leq |b_{i-1}|+|b_i|+|a_i|\,.
            \end{equation}
            The bound \eqref{eq:bornelambda} also holds for some $i$ for $|\lambda_\mathrm{min}|$, and because $\rho(J)=\max\{|\lambda_\mathrm{min}|,|\lambda_\mathrm{max}|\}$, we conclude on the second inequality.
		\end{Pro}

\begin{Rem}\label{rem:matricelower}		
 We observe that Lemma \ref{lem:inclusion_bourrine} allows us to conclude on the LD lower bound for $x\geq 1$ in Theorem \ref{thm:LDPtridiag_intro}. Indeed, by Lemma \ref{lem:inclusion_bourrine}, we have for $x\geq 1$
 \begin{equation}
     \label{eq:ldp_matrice_lower}
     \liminf_{N\rightarrow\infty}\frac{1}{\log N}\log\P\left(\max_{1\leq i \leq N} a_i\geq x\sqrt{2\log N}\right) \leq \liminf_{N\rightarrow\infty}\frac{1}{\log N}\log\P\left(\lambda_\mathrm{max}(J) \geq x\sqrt{2\log N}\right),
 \end{equation}
        and since the $a_i$'s are independent with uniform Gaussian tail, one establishes that the left hand side is $-I_2(x)$. Furthermore, note that for $x<1$,
\begin{equation}
     \label{eq:ldp_matrice_left}
     \limsup_{N\rightarrow\infty}\frac{1}{\log N}\log\P\left(\lambda_\mathrm{max}(J) \leq x\sqrt{2\log N}\right)\leq\limsup_{N\rightarrow\infty}\frac{1}{\log N}\log\P\left(\max_{1\leq i \leq N} a_i\leq x\sqrt{2\log N}\right).
 \end{equation}
 As before, the right hand side is $-\infty$, and we conclude on the left LD for the top eigenvalue.
 \end{Rem}

We will use the following result for the $L^2$ norm of independent variables with Gaussian tails.
\begin{Lem}[Gaussian tail estimate]\label{lem:gaussian_tail}
    Let $(X_i)_{i\geq 1}$ be a sequence of independent variables with uniform Gaussian tails as in Definition \ref{def:gaussian_tail}. Let $d\geq 1$ and $\mcal{I}_d$ be the set of subsets $I\subset\llbracket1,+\infty\llbracket$ with cardinal $\mrm{Card}(I)\leq d$. For all $I\in \mathcal{I}_d$, let
    $$Z_I=\left(\sum_{i\in I} X_i^2\right)^{1/2}$$
    be the euclidean norm of $(X_i)_{i\in I}$. Then, the family $(Z_I)_{I\in \mathcal{I}_d}$ has uniform Gaussian tail.
    \end{Lem}


\begin{Pro}[of Lemma \ref{lem:gaussian_tail}]
We start with the case $d=2$. Let $i_1,i_2\in I$ and $X:=X_{i_1}$, $Y:=X_{i_2}$. We want to show that the variable $Z:=(X^2+Y^2)^{1/2}$ satisfies as $x\rightarrow+\infty$:
$$ -\log \P(Z>x) = x^2/2 +o(x^2),$$
where the remainder is bounded independently of $(i_1,i_2)\in I^2$. Let $n\geq 1$ and $x>0$, we decompose
\begin{align*}
    \P(|X|>x)\leq\P(Z^2>x^2)&\leq\P(|X|>x)+ \sum_{k=0}^{n-1}\P\left(X^2\in \left[\frac{k}{n}x^2,\frac{k+1}{n}x^2\right], Y^2>x^2\left(1-\frac{k+1}{n}\right)\right)
    \\&\defi \P(|X|>x) + S_n.
\end{align*}
By assumption, $\log \P(|X|>x)=-x^2/2 + o(x^2),$ where the remainder is independent of $i_1\in I$. We now establish that $S_n\leq e^{-x^2/2 +o(x^2)}$ in the previous decomposition, with $o(x^2)$ bounded independently of $i_1,i_2$. By independence of $X,Y$, we have
\begin{multline*}
    S_n \leq \sum_{k=0}^{n-1}\P(X^2\geq kx^2/n, Y^2\geq x^2(1-(k+1)/n)) = \sum_{k=0}^{n-1}\P(X^2\geq kx^2/n) \P(Y^2\geq x^2(1-(k+1)/n))\\
    \leq n\max_{0\leq k \leq n-1} \P(X^2\geq kx^2/n)\P(Y^2\geq x^2(1-(k+1)/n)).
\end{multline*}
Set $n=\lfloor x \rfloor$. We have by assumption for $0\leq k \leq n-1=\lfloor x \rfloor -1$, and using that $1\leq x/\lfloor x \rfloor \leq 1+1/\lfloor x\rfloor$,
\begin{multline*}
    -\log\P(X^2\geq kx^2/n)-\log\P(Y^2\geq x^2(1-(k+1)/n)) = \frac{kx}{2} +o(kx) + \frac{x^2-x(k+1)(1+1/\lfloor x \rfloor)}{2} + o(x^2-kx) \\= \frac{x^2}{2} + o(x^2),
\end{multline*}
where we used that $k\leq \lfloor x \rfloor$, thus both terms $o(kx)$ and $o(x^2-kx)$ are bounded by some $o(x^2)$ independent of $k$. Thus, $-\log S_n\leq -\log \lfloor x \rfloor + x^2/2 + o(x^2)=x^2/2 + o(x^2)$.
Therefore, adding up both estimates we have
$$ \P(|X|>x) = e^{-x^2/2 + o(x^2)} \leq \P(Z>x) \leq e^{-x^2/2 +o(x^2)} + e^{-x^2/2 +o(x^2)}=e^{-x^2/2 + o(x^2)},$$
where the terms $o(x^2)$ are independent of the indices $i_1,i_2$. This establishes the case $d=2$. In the general case, $Z^2=(X_{i_1}^2+\ldots + X_{i_{d-1}}^2) + X_{i_d}^2$. By induction, $X\defi(X_{i_1}^2+\ldots + X_{i_{d-1}}^2)^{1/2}$ has Gaussian tail, uniformly in $(i_1,\ldots,i_{d-1})\in I^{d-1}$, and $Y\defi X_{i_d}^2$ has Gaussian tail uniformly in $i_d\in I$. We conclude by the case $d=2$.
\end{Pro}

For $x\geq 1$, the upper bound on $\lambda_\mathrm{max}$ in Lemma \ref{lem:inclusion_bourrine} is too crude to capture the LD upper bound. Indeed, if $X,Y$ are standard Gaussians, then $|X|+|Y|$ doesn't have Gaussian tail\footnote{Instead, by a Laplace principle, one can check that $\lim_{a\to +\infty} -\log\mathbb{P}(|X|+|Y|\geq a)/a^2=\inf_{|u|+|v|\geq 1}(u^2+v^2)=1/4$, attained at $|u|=|v|=1/2$. An interpretation of this fact is saying that the less improbable way to make $|X|+|Y|$ big, is splitting evenly the mass between $|X|$ and $|Y|$.}. Instead, we will show that we can assume, at the scale of LD, that only a few entries of $J$ matter. Then we will be able to give a precise LD upper bound, using Lemma \ref{lem:gaussian_tail}.\\ 

For $\varepsilon>0$, and $J=T$ or $J=T_\mathsf{per}$, we decompose $J = J^\varepsilon + S^\varepsilon$, where 
\begin{equation}
    \label{def:J_eps}
    J^\varepsilon \defi (J_{i,j}\mathbf{1}_{ \{|J_{i,j}|\geq \varepsilon\sqrt{2\log N}\} })_{1\leq i,j\leq N}
\end{equation}
is the matrix where we cut the entries smaller than $\varepsilon\sqrt{2\log N}$, and $$S^\varepsilon \defi J-J^\varepsilon = (J_{i,j}\mathbf{1}_{ \{|J_{i,j}|< \varepsilon\sqrt{2\log N}\} })_{1\leq i,j \leq N}.$$  We have the following exponential approximation for the top eigenvalue of $J=J_N$.

\begin{Lem}
    \label{lem:exp_approx}
    Let $\overline{\lambda}_{\max} \defi \lambda_\mathrm{max}(J_N)/\sqrt{2\log N}$, and $\overline{\lambda}_{\max}^{(\varepsilon)}\defi \lambda_\mathrm{max}(J^\varepsilon_N)/\sqrt{2\log N}$. Then,
    $$ |\overline{\lambda}_{\max}-\overline{\lambda}_{\max}^{(\varepsilon)}|\leq 3\varepsilon.$$
    As a consequence, we have the exponential approximation
    $$ \limsup_{N\rightarrow\infty}\frac{1}{\log N}\log\P\left( |\overline{\lambda}_{\max}-\overline{\lambda}_{\max}^{(\varepsilon)}|>\delta \right) \underset{\varepsilon\to 0}{\longrightarrow} -\infty.$$
\end{Lem}

\begin{Pro}
Using that for real symmetric matrices one has
$$ \lambda_\mathrm{max}(A)=\max_{\|x\|=1} \braket{Ax,x},$$
one sees that $ \left|\lambda_\mathrm{max}(J)-\lambda_\mathrm{max}(J^\varepsilon)\right|\leq \rho(S^\varepsilon).$ Now, by Lemma \ref{lem:inclusion_bourrine}, we have $$\rho(S^\varepsilon_N)\leq \max_{1\leq i \leq N}\{ |S^\varepsilon_{i,i}|+|S^\varepsilon_{i,i-1}|+|S^\varepsilon_{i-1,i}|\}\leq 3\varepsilon\sqrt{2\log N}.$$
Now, for any $\delta>0$, we deduce that for $\varepsilon<\delta/3$
    $$ \P\left( |\overline{\lambda}_{\max}-\overline{\lambda}^{(\varepsilon)}_{\max}|>\delta \right) = 0.$$\qedsymbol{}
\end{Pro}

Let us now note that for a tridiagonal matrix $J=T(\underline{a},\underline{b})$ having vanishing off-diagonal entries $b_{i_1},\ldots, b_{i_k}=0$, we can decompose $J$ into sub blocks. As an example, see Figure \ref{fig:blocks}.

\begin{figure}[ht]
  \centering
  \begin{tikzpicture}[
  every node/.style={anchor=center, font=\large, inner sep=1pt},
  matrix cell/.style={nodes={minimum height=0.8em, inner sep=1pt}},
]

  \matrix[matrix of math nodes,
          left delimiter={(}, right delimiter={)},
          row sep=4pt, column sep=4pt,
          matrix cell] (m) {
    a_1  & b_1  &     &      &      &      &      &      &      &      \\
    b_1  & a_2  & 0  &    &      &      &      &      &      &      \\
       & 0  & a_3  & 0  &     &      &      &      &      &      \\
         &     & 0  & a_4  & b_4  &     &      &      &      &      \\
         &      &     & b_4  & a_5  & b_5  &     &      &      &      \\
         &      &      &     & b_5  & a_6  & 0  &     &      &      \\
         &      &      &      &     & 0  & a_7  & b_7  &     &      \\
         &      &      &      &      &     & b_7  & a_8  & b_8  &     \\
         &      &      &      &      &      &     & b_8  & a_9  & b_9  \\
         &      &      &      &      &      &      &     & b_9  & a_{10} \\
  };

  \begin{scope}[on background layer]
    \node[draw=blue, thick, rounded corners, fit=(m-1-1)(m-2-2), inner sep=2pt] {};
    
    \node[draw=blue, thick, rounded corners, fit=(m-3-3), inner sep=4pt] {};

    \node[draw=blue, thick, rounded corners, fit=(m-4-4)(m-6-6), inner sep=2pt] {};

    \node[draw=blue, thick, rounded corners, fit=(m-7-7)(m-10-10), inner sep=4pt] {};
  \end{scope}

\end{tikzpicture}

  \caption{Example where $N=10$, $b_2=b_3=b_6=0$.}
  \label{fig:blocks}
\end{figure}


\begin{Def}[Block decomposition]
    Let $J=T(\underline{a},\underline{b})$ be a symmetric tridiagonal matrix of size $N\times N$. Assume that $J$ has exactly $k$ vanishing off-diagonal entries:
    $$ b_{i_1},b_{i_2},\ldots,b_{i_k} = 0\qquad \text{and}\qquad b_i  \neq 0 \text{ otherwise.} $$
    We write the block decomposition $ J=\mathrm{diag}(B^{(1)},\ldots,B^{(k+1)}),$
    where $B^{(\ell)}=(J_{i,j})_{i_{\ell-1}< i,j\leq i_{\ell}}$, where we set $i_0=1$, $i_{k+1}=N$.\\
    Furthermore, we denote $d_\ell = i_\ell-i_{\ell-1}$ the dimension of $B_\ell$, and 
    \begin{equation}
\label{def:d_max}
    d_\mathrm{max}(J) = \max_{1\leq \ell \leq k+1} d_\ell = \max_{1\leq \ell \leq k+1} \big\{i_\ell-i_{\ell-1}\big\}.
\end{equation}
\end{Def}

We will now show that in the nonperiodic case, in the regime of LD, we can assume that $d_\mathrm{max}(J_N^\varepsilon)$ is finite.



\begin{Lem}\label{lem:block}
    Under the hypothesis of Theorem \ref{thm:LDPtridiag_intro}, let $J=T(\underline{a}^{(N)},\underline{b}^{(N)})$ tridiagonal. Let $\varepsilon>0$, and let $J_N^\varepsilon$ given by \eqref{def:J_eps}. Then, under the Gaussian tail assumption of Theorem \ref{thm:LDPtridiag_intro},
    $$ \limsup_{N\rightarrow\infty}\frac{1}{\log N}\log \P\left( d_\mathrm{max}(J_N^\varepsilon)\geq d \right) \underset{d \to +\infty}{\longrightarrow}-\infty.$$
\end{Lem}

\begin{Pro}
Let $d\geq1$, the size of the blocks in the block decomposition of $J_N^\varepsilon$ is determined by the vanishing off-diagonal entries: on the event $\{d_\text{max}(J^\varepsilon_N)\geq d\}$, there is at least some $1\leq i\leq N-d+1$ such that $|b_i^{(N)}|,|b_{i+1}^{(N)}|,\ldots,|b_{i+d-1}^{(N)}|\geq \varepsilon\sqrt{2\log N}$. Thus by a union bound on the starting position $i$ of such a sequence, independence and uniform Gaussian tail, we can assume without loss of generality that $i=1$, and
\begin{align*}
    \log\P\left( d_\mathrm{max}(J_N^\varepsilon)\geq d\right)&\leq \log N + \sum_{k=1}^{d}\log \P\left(|b_{k}^{(N)}|>\varepsilon\sqrt{2\log N}\right)\\
    &\leq \log N + \sum_{k=1}^{d}\big(-\varepsilon^2\log N +o(\log N)\big)= \log N - d\varepsilon^2\log N + o(\log N).
\end{align*}
Dividing by $\log N$ and letting $d\to +\infty$, we deduce the claim.
\end{Pro}

We now have everything to prove Theorem \ref{thm:LDPtridiag_intro}.

\begin{Pro}[of Theorem \ref{thm:LDPtridiag_intro}]
We will write $a_i:=a_i^{(N)}$ and $b_i:=b_i^{(N)}$, making implicit the $N$-dependence of the entries.
Recall that the LD lower bound for $x\geq 1$ was already established in \eqref{eq:ldp_matrice_lower}. We now establish the LD upper bound for $x\geq 1$. \\
By the exponential approximation of Lemma \ref{lem:exp_approx}, recalling the notation $\overline\lambda_{\mathrm{max}}= \lambda_\mathrm{max}(J)/\sqrt{2\log N}$, we have
\begin{equation}\label{eq:approx_eps}
    \limsup_{N\rightarrow\infty}\frac{1}{\log N}\log \P\left( \overline{\lambda}_{\max}>x\right) \leq \lim_{\delta\to 0}\liminf_{\varepsilon\to 0}\limsup_{N\rightarrow\infty}\frac{1}{\log N}\log \P\left( \overline{\lambda}^{(\varepsilon)}_{\max} > x-\delta \right).
\end{equation}
We now conclude on the nonperiodic tridiagonal case $J=T(\underline{a},\underline{b})$. 

For fixed $\varepsilon>0$, we have for any $d\geq 1$, recalling the notation $d_\mathrm{max}$ \eqref{def:d_max},
\begin{multline*}
    \limsup_{N\rightarrow\infty}\frac{1}{\log N}\log \P\left( \overline{\lambda}_{\max}^{(\varepsilon)} > x \right)
    =\max\Big\{ \limsup_{N\rightarrow\infty}\frac{1}{\log N}\log\P\left( \overline{\lambda}_{\max}^{(\varepsilon)} > x, d_{\max}(J_N^\varepsilon)\leq d\right)\ ;\\ \limsup_{N\rightarrow\infty}\frac{1}{\log N}\log \P\left( \overline{\lambda}_{\max}^{(\varepsilon)} > x, d_{\max}(J_N^\varepsilon)> d\right)\Big\}.
\end{multline*}
Using Lemma \ref{lem:block}, taking $d\to +\infty$,
\begin{equation}
    \label{eq:approx_d}
    \limsup_{N\rightarrow\infty}\frac{1}{\log N}\log \P\left( \overline{\lambda}_{\max}^{(\varepsilon)} > x \right) = \limsup_{d\to +\infty}\limsup_{N\rightarrow\infty}\frac{1}{\log N}\log\P\left( \overline{\lambda}_{\max}^{(\varepsilon)} > x, d_{\max}(J_N^\varepsilon)\leq d\right).
\end{equation}
Now, if $\lambda\in \mathrm{Sp}(J^\varepsilon)$, there is some $1\leq \ell \leq k_N^\varepsilon$ such that $\lambda \in \mathrm{Sp}(B_\ell)$. Then, on $\left\{ d_\mathrm{max}(J^\varepsilon)\leq d\right\},$
\begin{align*}
    |\lambda| \leq \tr(B_\ell^2)^{1/2} &= \Big(\sum_{j=1}^{d_\ell} B_\ell(j,j)^2+2B_\ell(j,j+1)^2\Big)^{1/2}\\ &= \Big(\sum_{j=1}^{d_\ell}J^{\varepsilon}(i_{\ell-1}+j,i_{\ell-1}+j)^2+2J^{\varepsilon}(i_{\ell-1}+j,i_{\ell-1}+j+1)^2\Big)^{1/2}\\
    &\leq\Big(\sum_{j=1}^d a_{i_{\ell-1}+j}^2+2b_{i_{\ell-1}+j}^2\Big)^{1/2}
    \defi Z_{i_\ell}^{(d)},
\end{align*}
where we used that $d_\ell\leq d$, that $|J^\varepsilon_N(i,j)|\leq |J_N(i,j)|$ and where for any index $1\leq i \leq N$ we defined
$$ Z_i^{(d)}=\Big(\sum_{j=1}^{d}a_{i+j}^2+b_{i+j}^2\Big)^{1/2},$$
where the indices are taken modulo $N$.
Thus $$ \P(\overline{\lambda}_{\max}^{(\varepsilon)}>x,d_{\max}(J_N^\varepsilon)\leq d)\leq \P\left(\max_{1\leq \ell \leq k_N^\varepsilon} Z_{i_\ell}^{(d)}>x\sqrt{2\log N}\right)\leq \P\left(\max_{1\leq i \leq N} Z_i^{(d)}>x\sqrt{2\log N}\right),$$
where we used $k_N^\varepsilon\leq N$.\\
Using that each $Z_i^{(d)}$ has Gaussian tail independent of $N$ by the fact that each $a_i$ and $\sqrt{2}b_i$ have Gaussian tail and Lemma \ref{lem:gaussian_tail}, we conclude by an union bound that for any fixed $d\geq 1$, and for $x\geq 1$,
\begin{align*}
    \limsup_{N\rightarrow\infty}\frac{1}{\log N}\log\P\left( \overline{\lambda}_{\max}^{(\varepsilon)} >x, d_\mathrm{max}\leq d\right) &\leq  -I_2(x).
\end{align*}
By equations \eqref{eq:approx_d} and \eqref{eq:approx_eps}, we conclude for the LD upper bound in the case $x\geq 1$. Finally, the left deviations (\textit{i.e.} the case $x<1$) were established by \eqref{eq:ldp_matrice_left}, which allows us to conclude on the full LDP in the nonperiodic case.

In the periodic case $J=T_\mathsf{per}(\underline{a},\underline{b})$, the only difference with the previous proof is when decomposing $J^{\varepsilon}$ into blocks. We will show that at the scale of large deviations at least one off diagonal entry of $J^\varepsilon$ vanishes, and perform a change of basis to deal with a tridiagonal matrix. Note that for any $\varepsilon>0$, 
$$\log\P(\forall i,\ |b_i|\geq \varepsilon\sqrt{2\log N})=\sum_{i=1}^N \log\P(|b_i|\geq \varepsilon\sqrt{2\log N})\underset{N\to +\infty}{\sim} -N\varepsilon^2\log N,$$
where we used the Gaussian tail of the entries. Thus, we see that
$$ \frac{1}{\log N}\log\P(\forall i,|b_i|\geq \varepsilon\sqrt{2\log N}) \underset{N\to +\infty}{\longrightarrow} -\infty.$$
As a consequence, by \eqref{eq:approx_eps},
$$ \limsup_{N\to \infty} \frac{1}{\log N}\log\P(\overline{\lambda}_{\max}>x) \leq \lim_{\delta\to 0}\liminf_{\varepsilon\to 0}\limsup_{N\to \infty}\frac{1}{\log N}\log \P(\overline{\lambda}_{\max}^{(\varepsilon)}>x-\delta, \mathcal{B}_N^{\varepsilon}),$$
where $\mathcal{B}_N^{\varepsilon}=\{\exists i\in \llbracket 1,N\rrbracket:|b_i|\leq\varepsilon\sqrt{2\log N} \}.$ Now, note that on $\mathcal{B}_N^{\varepsilon}$, there exists $i^*$ such that $J^\varepsilon(i^*,i^*+1)=0$, therefore denoting $(e_1,\ldots,e_N)$ the canonical basis of $\R^N$, representing $J^\varepsilon$ in the shifted (random) basis $(e_{i^*+1},e_{i^*+2},\ldots,e_N,e_1,\ldots,e_{i^*})$, we have
\begin{equation}\label{eq:decompositionorthogonale}
    J^\varepsilon=UT^\varepsilon(\underline{a}^\prime,\underline{b}^\prime)U^\top,
\end{equation}
where $U$ is the change of basis matrix and $\underline{a}^\prime$, $\underline{b}^\prime$ are obtained from $\underline{a}$ and $\underline{b}$ by shifting the coordinates by $i^*$ modulo $N$ (see Figure \ref{fig2}). 

Coming back to the general case, $T^\varepsilon:=T^\varepsilon(\underline{a}^\prime,\underline{b}^\prime)$ and $J^\varepsilon=J(\underline{a}^\prime,\underline{b}^\prime)$ have the same spectrum, thus the same top eigenvalue. Now by the previous reasoning, we find the same bound $\lambda_{\max}(T^\varepsilon)=\lambda_{\max}(J^\varepsilon)\leq \max_{1\leq i \leq N} Z_i^{(d)}$, and we conclude as in the nonperiodic case.
\end{Pro}

\begin{figure}[ht]
  \centering
    $$ 
    \begin{pmatrix}
        a_1  & b_1  &     &      & b_5\\
    b_1  & a_2  & b_2 &      &    \\
         & b_2  & a_3 & 0  &    \\
         &      & 0 & a_4  & b_4\\
    b_5  &      &     & b_4  & a_5
    \end{pmatrix}=
    \begin{pmatrix}
          &   & 1    &      & \\
          &   &      & 1    &    \\
         &   &  &   &  1  \\
        1 &      &  &   & \\
      &  1    &     &   & 
    \end{pmatrix}
    \begin{pmatrix}
        a_4  & b_4  &     &      & \\
    b_1  & a_5  & b_5 &      &    \\
         & b_5  & a_1 & b_1  &    \\
         &      & b_1 & a_2  & b_2\\
      &      &     & b_2  & a_3
    \end{pmatrix} \begin{pmatrix}
          &   & 1    &      & \\
          &   &      & 1    &    \\
         &   &  &   &  1  \\
        1 &      &  &   & \\
      &  1    &     &   & 
    \end{pmatrix}^\top
    $$
    \caption{Illustration of the change of basis \eqref{eq:decompositionorthogonale} in the case $N=5$. In this example, $b_3=0$.}
    \label{fig2}
\end{figure}

	
	\bibliographystyle{alpha}
	
	\bibliography{BiblioLargest}

@article{pakzad2020large,
  title={Large Deviations Principle for the Largest Eigenvalue of the {G}aussian $\beta$-Ensemble at High Temperature},
  author={Pakzad, C.},
  journal={Journal of Theoretical Probability},
  year={2020}
}

@article{BordCaputo,
author = {Bordenave, C. and Caputo, P.},
title = {{A large deviation principle for {W}igner matrices without {G}aussian tails}},
journal = {The Annals of Probability},
year = {2014}
}

@article{Augeri,
author = {Augeri, F.},
title = {{Large deviations principle for the largest eigenvalue of {W}igner matrices without {G}aussian tails}},
volume = {21},
journal = {Electronic Journal of Probability},
year = {2016}
}

@article{vivo2015large,
  title={Large deviations of the maximum of independent and identically distributed random variables},
  author={Vivo, P.},
  journal={European Journal of Physics},
  year={2015}
}

@article{lambert2021poisson,
  title={Poisson statistics for {G}ibbs measures at high temperature},
  author={Lambert, G.},
  year={2021},
  journal={Ann. Inst. Henri Poincar\'{e}}
}

@article{basak2023upper,
  title={Upper tail of the spectral radius of sparse {E}rd{\"o}s--{R}{\'e}nyi graphs},
  author={Basak, A.},
  journal={Probability Theory and Related Fields},
  year={2023},
  publisher={Springer}
}

@article{cook2020large,
  title={Large deviations of subgraph counts for sparse {E}rd{\H{o}}s--{R}{\'e}nyi graphs},
  author={Cook, N. and Dembo, A.},
  journal={Advances in Mathematics},
  year={2020},
  publisher={Elsevier}
}

@article{cook2023full,
  title={Full large deviation principles for the largest eigenvalue of sub-Gaussian Wigner matrices},
  author={Cook, N.A. and Ducatez, R. and Guionnet, A.},
  journal={arXiv preprint arXiv:2302.14823},
  year={2023}
}

@book{dembo2009large,
  title={Large deviations techniques and applications},
  author={Dembo, A. and Zeitouni, O.},
  year={2009},
  publisher={Springer}
}

@article{dean2008extreme,
	title={Extreme value statistics of eigenvalues of {G}aussian random matrices},
	author={Dean, D. and Majumdar, S.},
	journal={Physical Review E},
	year={2008}
}

@article{ganguly2024spectral,
  title={Spectral large deviations of sparse random matrices},
  author={Ganguly, S. and Hiesmayr, E. and Nam, K.},
  journal={Journal of the London Mathematical Society},
  year={2024},
}

@article{Dwo,
	title={Asymptotics of the partition function for $\beta$-ensembles at high temperature},
	author={Dworaczek Guera, C.},
	volume = {31},
	journal = {Electronic Journal of Probability},
    year={2026}
}

@article {NakanoTrinh,
	AUTHOR = {Nakano, F. and Trinh, K. D.},
	TITLE = {{G}aussian $\beta$-ensembles at high temperature: eigenvalue fluctuations and bulk statistics.},
	JOURNAL = {J.Stat.Phys},
	YEAR = {2018}
}

@article{benaych2015poisson,
	title={{P}oisson statistics for matrix ensembles at large temperature},
	author={Benaych-Georges, F. and P{\'e}ch{\'e}, S.},
	journal={Journal of Statistical Physics},
	year={2015}
}

@article{pakzad2018poisson,
	title={{P}oisson statistics at the edge of {G}aussian $\beta$-ensembles at high temperature},
	author={Pakzad, C.},
	journal={arXiv preprint arXiv:1804.08214},
	year={2018}
}

@article{AllezBouchaudGuionnet,
	KEY={ABG12},
	AUTHOR = {Allez, R. and Bouchaud, J.P. and Guionnet, A.},
	TITLE = {Invariant $\beta$-ensembles and the {G}auss-{W}igner Crossover},
	JOURNAL = {Physical Review Letters},
	DOI = {109(9):094102},
	YEAR = {2012},
}

@Article{NakanoTrinh20,
	Author = {Nakano, F. and Trinh, K. D.},
	Title = {{P}oisson statistics for $\beta$-ensembles on the real line at high temperature},
	Journal = {J. Stat. Phys.},
	Year = {2020}
}

@article{bhattacharya2020upper,
  title={Upper tails for edge eigenvalues of random graphs},
  author={Bhattacharya, B. B. and Ganguly, S.},
  journal={SIAM Journal on Discrete Mathematics},
  year={2020}
}

@article{fan2015convergence,
  title={Convergence of eigenvalues to the support of the limiting measure in critical $\beta$ matrix models},
  author={Fan, C. and Guionnet, A. and Song, Y. and Wang, A.},
  journal={Random Matrices: Theory and Applications},
  year={2015}}

@article{allez2014tracy,
  title={Tracy--{W}idom at high temperature},
  author={Allez, R. and Dumaz, L.},
  journal={Journal of Statistical Physics},
  year={2014},
  publisher={Springer}
}

@article{allez2014sine,
  title={From {S}ine kernel to {P}oisson statistics},
  author={Allez, R. and Dumaz, L.},
  year={2014},
  journal={Electronic Journal of Probability}
}

@article{johansson2007gumbel,
  title={From {G}umbel to {T}racy-{W}idom},
  author={Johansson, K.},
  journal={Probability theory and related fields},
  year={2007}
}

@article{borot2013asymptotic,
  title={Asymptotic expansion of $\beta$ matrix models in the one-cut regime},
  author={Borot, G. and Guionnet, A.},
  journal={Communications in Mathematical Physics},
  year={2013}
}

@article{guionnet2020large,
  title={Large deviations for the largest eigenvalue of {R}ademacher matrices},
  author={Guionnet, A. and Husson, J.},
  year={2020},
  journal={The Annals of Probability}
}

@article{husson2022large,
  title={Large deviations for the largest eigenvalue of matrices with variance profiles},
  author={Husson, J.},
  journal={Electronic Journal of Probability},
  year={2022},
}

@article{nakano2024classical,
  title={Classical $\beta$-ensembles and related eigenvalues processes at high temperature and the {M}arkov-{K}rein transform},
  author={Nakano, F. and Trinh, H-D. and Trinh, K-D.},
  journal={arXiv preprint arXiv:2412.01015},
  year={2024}
}

@article{askey1984associated,
  title={Associated {L}aguerre and {H}ermite polynomials},
  author={Askey, R. and Wimp, J.},
  journal={Proceedings of the Royal Society of Edinburgh Section A: Mathematics},
  year={1984}
}

@article{ducatez2024large,
  title={Large deviation principle for the largest eigenvalue of random matrices with a variance profile},
  author={Ducatez, R. and Guionnet, A. and Husson, J.},
  journal={arXiv preprint arXiv:2403.05413},
  year={2024}
}

@article{hardy2021clt,
	title={{CLT} for circular $\beta$-ensembles at high temperature},
	author={Hardy, A. and Lambert, G.},
	journal={Journal of Functional Analysis},
	year={2021}
}

@inproceedings{mazzuca2024clt,
	title={{CLT} for $\beta$-Ensembles at High Temperature and for Integrable Systems: A Transfer Operator Approach},
	author={Mazzuca, G. and Memin, R.},
	booktitle={Ann. Inst. Henri Poincar\'{e}},
	year={2024},
}

@article {GMToda,
	AUTHOR = {Guionnet, A. and Memin, R.},
	TITLE = {Large deviations for {G}ibbs ensembles of the classical {T}oda chain},
	JOURNAL = {Electron. J. Probab.},
	YEAR = {2022}
}

@article{DwoMemin,
	author = {Dworaczek Guera,  C. and Memin, R.},
	title = {{CLT for real $\beta$-ensembles at high temperature}},
	journal = {Electronic Journal of Probability},
	year = {2024}}

@book{anderson2010introduction,
	title={An introduction to random matrices},
	author={Anderson, G. W. and Guionnet, A. and Zeitouni, O.},
	year={2010},
	publisher={Cambridge university press},
}

@article{padilla2025poisson,
  title={Poisson Statistics for {C}oulomb Gases at Intermediate Temperature Regimes},
  author={Padilla-Garza, D. and Peilen, L. and Thoma, E.},
  journal={arXiv preprint arXiv:2507.08198},
  year={2025}
}

@article {Garcia,
	AUTHOR = {Garc\'{\i}a-Zelada, D.},
	TITLE = {A large deviation principle for empirical measures on {P}olish
	spaces: application to singular {G}ibbs measures on manifolds},
	JOURNAL = {Ann. Inst. Henri Poincar\'{e}},
	YEAR = {2019}
}

@article{guionnet2024asymptoticexpansionpartitionfunction,
      title={Asymptotic expansion of the partition function for $\beta$-ensembles with complex potentials}, 
      author={A. Guionnet and K. Kozlowski and A. Little},
      year={2024},
      journal={arXiv preprint arXiv:2411.10610}
}

@article{maida2007large,
  title={Large deviations for the largest eigenvalue of rank one deformations of {G}aussian ensembles},
  author={Ma{\"\i}da, M.},
journal={Electronic Journal of Probability},
  year={2007}
}

@article{benaych2012large,
  title={Large deviations of the extreme eigenvalues of random deformations of matrices},
  author={Benaych-Georges, F. and Guionnet, A. and Ma{\"\i}da, M.},
  journal={Probability Theory and Related Fields},
  year={2012}
}

@article{arous1997large,
  title={Large deviations for {W}igner's law and {V}oiculescu's non-commutative entropy},
  author={Ben-Arous, G. and Guionnet, A.},
  journal={Probability theory and related fields},
  year={1997}
}

@article{arous2001aging,
  title={Aging of spherical spin glasses},
  author={Ben-Arous, G.  and Dembo, A. and Guionnet, A.},
  journal={Probability theory and related fields},
  year={2001}
}

@article{augeri2021large,
  title={Large deviations for the largest eigenvalue of sub-{G}aussian matrices},
  author={Augeri, F. and Guionnet, A. and Husson, J.},
  journal={Communications in mathematical physics},
  year={2021}
}

@article{ma2023unified,
  title={Unified limits and large deviation principles for $\beta$-{L}aguerre ensembles in global regime},
  author={Ma, Y-T.},
  journal={Acta Mathematica Sinica, English Series},
  year={2023}
}

@article{lei2023large,
  title={Large deviations for top eigenvalues of $\beta$-{J}acobi ensembles at scaling temperatures},
  author={Lei, L. and Ma, Y-T.},
  journal={Acta Mathematica Scientia},
  year={2023}
}

@article{padilla2024emergence,
  title={Emergence of a {P}oisson process in weakly interacting particle systems},
  author={Padilla-Garza, D. and Peilen, L. and Thoma, E.},
  journal={arXiv preprint arXiv:2405.02625},
  year={2024}
}

@article{ma2023limit,
  title={Limit Theorems and Large Deviations for $\beta$-{J}acobi Ensembles at Scaling Temperatures},
  author={Ma, Y-T.},
  journal={Acta Mathematica Sinica, English Series},
  year={2023}
}

@article{ganguly2022large,
  title={Large deviations for the largest eigenvalue of {G}aussian networks with constant average degree},
  author={Ganguly, S. and Nam, K.},
  journal={Probability Theory and Related Fields},
  year={2022}
}

@article{bhattacharya2021spectral,
  title={Spectral edge in sparse random graphs: Upper and lower tail large deviations},
  author={Bhattacharya, B. B. and Bhattacharya, S. and Ganguly, S.},
  year={2021},
  journal={The Annals of Probability}
}

@article {MMS,
    AUTHOR = {Ma{\"\i}da, M. and Maurel-Segala, E.},
     TITLE = {Free transport-entropy inequalities for non-convex potentials
              and application to concentration for random matrices},
   JOURNAL = {Probab. Theory Related Fields},
      YEAR = {2014},
 
}

@article{benaych2019largest,
  title={Largest eigenvalues of sparse inhomogeneous {E}rd{\H{o}}s--{R\'e}nyi graphs},
  author={Benaych-Georges, F. and Bordenave, C. and Knowles, A.},
  journal={The Annals of Probability},
  year={2019},
}

@article{benaych2020spectral,
  title={Spectral radii of sparse random matrices},
  author={Benaych-Georges, F. and Bordenave, C. and Knowles, A.},
JOURNAL = {Ann. Inst. Henri Poincar\'{e} Probab. Stat.},
  year={2020}
}

@article{alt2023poisson,
  title={Poisson statistics and localization at the spectral edge of sparse {E}rd{\H{o}}s--{R\'e}nyi graphs},
  author={Alt, J. and Ducatez, R. and Knowles, A.},
  journal={The Annals of Probability},
  year={2023},
}

@article{augeri2023large,
  title={Large deviations of the largest eigenvalue of supercritical sparse {W}igner matrices},
  author={Augeri, F. and Basak, A.},
  journal={The Annals of Probability},
  year={2023}
}

@article{spohn2020generalized,
  title={Generalized {G}ibbs ensembles of the classical {T}oda chain},
  author={Spohn, H.},
  journal={Journal of Statistical Physics},
  year={2020},
}

@article{altducatezknowles,
author = {Alt, J. and Ducatez, R. and Knowles, A.},
title = {{Extremal eigenvalues of critical {E}rdős–{R}ényi graphs}},
journal = {The Annals of Probability},
year = {2021},

}

@article{dued,
    AUTHOR = {Dumitriu, I. and Edelman, A.},
     TITLE = {Matrix models for $\beta$-ensembles},
   JOURNAL = {J. Math. Phys.},
      YEAR = {2002},
}

@article{mazzuca2022mean,
  title={On the mean density of states of some matrices related to the $\beta$-ensembles and an application to the {T}oda lattice},
  author={Mazzuca, G.},
  journal={Journal of Mathematical Physics},
  year={2022}
}

@article{kerov1998interlacing,
  title={Interlacing measures},
  author={Kerov, S.},
  journal={American Mathematical Society Translations},
  pages={35--84},
  year={1998},
  publisher={American Mathematical Society}
}

@article{grava2026large,
  title={Large deviations of the periodic Toda chain},
  author={Grava, T. and Guionnet, A. and Kozlowski, K.K. and Little, A.},
  journal={arXiv preprint arXiv:2604.00635},
  year={2026}
}

@inproceedings{han2025deviation,
  title={Deviation of Top Eigenvalue for Some Tridiagonal Matrices Under Various Moment Assumptions},
  author={Han, Y.},
  booktitle={Annales Henri Poincar{\'e}},
  volume={26},
  number={6},
  pages={1907--1926},
  year={2025},
  organization={Springer}
}

@article{augeri2025large,
  title={Large deviations of the empirical spectral measure of supercritical sparse Wigner matrices},
  author={Augeri, F.},
  journal={Advances in Mathematics},
  volume={466},
  pages={110156},
  year={2025},
  publisher={Elsevier}
}
	
\end{document}